\documentclass[12pt]{article}
\usepackage{amssymb,amsmath,amsfonts}
\usepackage{relsize}
\usepackage[sort]{cite}
\usepackage{esvect}
\usepackage{txfonts}
\usepackage{graphicx}

\newtheorem{theorem}{Theorem}[section]
\newtheorem{definition}[theorem]{Definition}

\newtheorem{lemma}[theorem]{Lemma}

\newtheorem{problem}[theorem]{Problem}

\newtheorem{assumption}[theorem]{Assumption}
\newtheorem{statement}[theorem]{Statement}
\newtheorem{convention}[theorem]{Convention}
\newcommand {\Ac}      {{\mathcal A}}

\newcommand {\Ec}      {{\mathcal E}}
\newcommand {\Fc}      {{\mathcal F}}
\newcommand {\Gc}      {{\mathcal G}}

\newcommand {\Ic}      {{\mathcal I}}
\newcommand {\Jc}      {{\mathcal J}}
\newcommand {\Kc}      {{\mathcal K}}
\newcommand {\Lc}      {{\mathcal L}}
\newcommand {\Mc}      {{\mathcal M}}

\newcommand {\Pc}      {{\mathcal P}}

\newcommand {\Tc}      {{\mathcal T}}

\newcommand {\Wc}      {{\mathcal W}}
\newcommand {\R}       {{\mathbb R}}

\newcommand {\N}       {{\mathbb N}}
\newcommand {\Z}       {{\mathbb Z}}

\newcommand {\tE}      {\widetilde{E}}
\newcommand {\tF}      {\widetilde{F}}

\newcommand {\tS}      {\widetilde{S}}

\newcommand {\tY}      {\widetilde{Y}}

\newcommand {\tf}      {\tilde{f}}

\newcommand {\tn}      {\tilde{n}}

\newcommand {\RN}      {\R^n}
\newcommand {\ve}      {\varepsilon}

\newcommand {\LMPR}    {L_{p}^{m}(\R)}
\newcommand {\WORN}    {W_{p}^{1}(\RN)}
\newcommand {\WMP}     {W_{p}^{m}(\RN)}
\newcommand {\WMPR}    {W_{p}^{m}(\R)}
\newcommand {\LMIR}    {L^m_\infty(\R)}

\newcommand {\WMIR}    {W^m_\infty(\R)}

\newcommand {\LPR}     {L_p(\R)}
\newcommand {\LIR}     {L_\infty(\R)}

\newcommand {\intl}    {\int\limits}
\newcommand {\emp}     {\emptyset}
\newcommand {\PM}      {\Pc_{m}}

\newcommand {\NMP}     {\Lc_{m,p}}
\newcommand {\TNMP}    {\widetilde{\Lc}_{m,p}}
\newcommand {\TNMI}    {\widetilde{\Lc}_{m,\infty}}

\newcommand {\NWMP}    {\Wc_{m,p}}
\newcommand {\NWMI}    {\Wc_{m,\infty}}
\newcommand {\TLNW}    {\widetilde{\Wc}_{m,p}}
\newcommand {\TLNWI}   {\widetilde{\Wc}_{m,\infty}}
\newcommand {\WCP}     {\widehat{\Wc}_{m,p}}

\newcommand {\NIN}     {\Lc}
\newcommand {\SH}      {{S\hspace*{-0.5mm}}}

\newcommand {\ME}      {m_E}
\newcommand {\meh}     {\hspace{0.1mm}}

\newcommand {\xb}      {\bar{x}}
\newcommand {\tlm}     {\widetilde{\lambda}}

\newcommand {\WOP}     {\Fc_{m,E}^{(Wh)}}
\newcommand {\WOPTE}   {\Fc_{m,\tE}^{(Wh)}}

\newcommand {\vkp}     {\varkappa}
\newcommand {\brf}     {\bar{f}}
\newcommand {\fks}     {f^{\#}_{k,E}}
\newcommand {\fjs}     {f^{\#}_{j,E}}
\newcommand {\fms}     {f^{\#}_{m,E}}

\newcommand {\smsk}    {\smallskip}
\newcommand {\msk}     {\medskip}
\newcommand {\bsk}     {\bigskip}

\newcommand {\capsm}   {\mathsmaller{\bigcap}}

\newcommand {\mcup}    {\mathlarger{\cup}}
\newcommand {\mcap}    {\mathlarger{\cap}}

\newcommand {\smed}    {\mathlarger{\sum}}
\newcommand {\diam}    {\operatorname{diam}}
\newcommand {\dist}    {\operatorname{dist}}

\newcommand {\EXT}     {\operatorname{Ext}}

\newcommand {\VST}     {\vspace*{1mm}}
\newcommand {\bx}      {\hspace{10mm}$\Box$}

\newcommand {\nn}      {\nonumber}
\newcommand {\xlabel}      {\label}

\newcommand {\rf}[1]    {(\ref{#1})}      
\newcommand {\reff}[1] {\ref{#1}}         
\newcommand{\lbl}[1]      {\label{#1}}       
\newcommand{\be}          {\begin{eqnarray}}
\newcommand{\bel}[1]      {\begin{eqnarray} \label{#1}}
\newcommand{\ee}           {\end{eqnarray}}
\newcommand {\SECT}[2] {\section*{\centerline{\normalsize
{\bf #1}}} \setcounter{section}{#2}
\setcounter{theorem}{0}\setcounter{equation}{0}}
\begin{document}
\parindent 1em
\parskip 0mm
\medskip
\centerline{\LARGE Sobolev functions on closed subsets of the real line}\vspace*{10mm}
\centerline{\large Pavel Shvartsman}\vspace*{5 mm}
\centerline {\it Department of Mathematics, Technion - Israel Institute of Technology, 32000 Haifa, Israel}\vspace*{2 mm}
\vspace*{20 mm}
\renewcommand{\thefootnote}{ }
\footnotetext[1]{\hspace{-6mm}
{\it E-mail address:} pshv@technion.ac.il}
\hrule\msk
\par\noindent {\bf Abstract}
\bsk
\par For each $p>1$ and each positive integer
$m$ we use divided differences to give intrinsic characte\-rizations of the restriction of the Sobolev space $\WMPR$ to an arbitrary closed subset of the real line.
\bsk
{\small
\smsk
}
\msk
\hrule
\renewcommand{\contentsname}{ }
\tableofcontents
\addtocontents{toc}{{\bf Contents}\vspace*{5mm}\par}
\SECT{1. Introduction.}{1}
\addtocontents{toc}{1. Introduction.\hfill \thepage\par\VST}

\indent\par In this paper we characterize the restrictions of Sobolev functions of one variable to an arbitrary closed subset of the real line. For each $m\in\N$ and each $p\in (1,\infty]$, we consider $\LMPR$, the standard homogeneous Sobolev space on $\R$. We identify $\LMPR$ with the space of all real valued functions $F$ on $\R$ such that the $(m-1)$-th derivative $F^{(m-1)}$ is absolutely continuous on $\R$ and the weak $m$th derivative $F^{(m)}\in L_p(\R)$. $\LMPR$ is seminormed by
$\| F\| _{\LMPR}= \| F^{(m)}\| _{L_p(\R)}$.
\smsk

As usual, we let $\WMPR$ denote the corresponding Sobolev space of all functions $F\in \LMPR$ whose derivatives on $\R$ of \textit{all orders up to $m$} belong to $\LPR$. This space is normed by
\[
\| F\| _{\WMPR}=\smed_{k=0}^m\,\, \| F^{(k)}\| _{\LPR}.
\]

 In this paper, we study the following extension problem.
\begin{problem}
\lbl{PR-MAIN} Let $p\in(1,\infty]$, $m\in\N$, and let $E$ be a closed subset of $\R$.
Let $f$ be a function on $E$.
\smsk

How can we decide whether $f$ extends to a function $F\in \WMPR$?

If such an $F$ exists, then how small can its Sobolev norm $\| F\| _{\WMPR}$ be?
\end{problem}

We denote the infimum of all these norms by $\| f\| _{\WMPR\vert  _E}$; thus
 \begin{equation}\xlabel{N-WMPR}
\| f\| _{\WMPR\vert  _E}
=\inf \{\| F\| _{\WMPR}:F\in\WMPR, F\vert  _{E}=f\}.
\end{equation}
 We refer to $\| f\| _{\WMPR\vert  _E}$ as the \textit{trace norm on $E$ of the function $f$} in $\WMPR$. This quantity provides the standard quotient space norm on \textit{the trace space} $\WMPR\vert  _{E}$ of all restrictions of $\WMPR$-functions to $E$, i.e.,~on the space
\[
\WMPR\vert  _{E}=\{f:E\to\R:\text{there exists}~~F\in
\WMPR\ \ \text{such that}\ \ F\vert  _{E}=f\}.
\]

The trace space $\LMPR\vert  _{E}$ is defined in an analogous way.
\smallskip

This paper can be considered as a sequel to our earlier
paper~\cite{Sh-LMP-2018} in which we dealt with the analo\-gous
problem for homogeneous Sobolev spaces.  The reader may notice
various analogies between material presented here and
in~\cite{Sh-LMP-2018}.

Theorem \reff{W-VAR-IN}, our main contribution in this paper, provides a complete solution to Problem \reff{PR-MAIN}. Let us prepare the ingredients that are needed to formulate this theorem: Given a function $f$ defined on a $k+1$ point set $S=\{x_0, \ldots ,x_k\}\subset\R$, we let $\Delta^kf[x_0, \ldots ,x_{k}]$ denote the $k^{\,\text{th}}$ \textit{order divided difference} of $f$ on $S$. (Recall that $\Delta^kf[x_0, \ldots ,x_{k}]$ coincides with the coefficient of $x^k$ in the Lagrange polynomial of degree at most $k$ which agrees with $f$ on $S$. See Section 2 for other equivalent definitions of divided differences and their main properties.)

Everywhere in this paper we will use the following notation: Given a finite strictly increasing sequence $\{x_0, \ldots ,x_n\}\subset E$ we set
 \begin{equation}\xlabel{AGR}
x_{i}=+\infty~~~\text{if}~~~i>n.
\end{equation}

Given $p\in(1,\infty)$, a set $E\subset\R$ with $\#E\ge m+1$, and a function $f$ on $E$, we set
 \begin{equation}\xlabel{N-WRS-IN}
\NWMP(f:E)=\sup_{\{x_0, \ldots ,x_n\}\subset E}
\,\,\left(\,\smed_{k=0}^{m}\,\,\smed_{i=0}^{n-k}
\min\left\{1,x_{i+m}-x_{i}\right\}
\left|\Delta^kf[x_i, \ldots ,x_{i+k}]\right|^p
\right)^{\frac1p}.
\end{equation}
 Here the supremum is taken over all integers $n\ge m$ and all strictly increasing sequences $\{x_0, \ldots ,x_n\}\subset E$ of $n+1$ elements.

We also put
 \begin{equation}\xlabel{N-WINF}
\NWMI(f:E)=\,\sup_{{k=0, \ldots ,m}\atop
{\{y_0, \ldots ,y_k\}\subset E}}
\,\,\,
\left|\Delta^kf[y_0, \ldots ,y_k]\right|
\end{equation}
 where the supremum is taken over all $k=0, ... ,m$, and all
strictly increasing sequences $\{y_0, ... ,y_k\}\subset E$.

Here now is the main result of our paper:
\begin{theorem}
\lbl{W-VAR-IN} Let $m$ be a positive integer, $p\in(1,\infty]$, and let $E$ be a closed subset of $\R$ containing at least $m+1$ points.

A function $f:E\to\R$ can be extended to a function $F\in\WMPR$ if and only if $\NWMP(f:E)<\infty$. Furthermore,
 \begin{equation}\xlabel{NM-CLC-IN}
\| f\| _{\WMPR\vert  _E}\sim \NWMP(f:E)\,,
\end{equation}
 where the constants in equivalence \rf{NM-CLC-IN} depend only on $m$.
\end{theorem}

We refer to Theorem \reff{W-VAR-IN} as a \textit{variational
extension criterion for the traces of $\WMPR$-functions}.
This theorem is proven in Section 3 and 4 of the present paper.
For variational extension criteria for the space $\WORN$,
$p>n$, and jet-spaces generated by functions from
$\WMP$, $p>n$, we refer the reader to~\cite[Theorem
1.2]{Sh2} and~\cite[Theorem 8.1]{Sh5} respectively.
\smallskip

An examination of our proof of Theorem \reff{W-VAR-IN} shows that when $E$ is a \textit{strictly increasing sequence} of points in $\R$ (finite, one-sided infinite, or bi-infinite) it is enough to take the supremum in \rf{N-WRS-IN} over a unique subsequence of $E$ -- the sequence $E$ itself. Therefore, in the particular case
of such sets $E$, we can obtain the following refinement of Theorem \reff{W-VAR-IN}:
\begin{theorem}
\lbl{W-TFIN} Let $\ell_1,\ell_2\in\Z\cup\{\pm\infty\}$, $\ell_1\le\ell_2$, and let $p\in(1,\infty)$. Let  $E=\{x_i\}_{i=\ell_1}^{\ell_2}$ be a strictly increasing sequence of points in $\R$, and let
 \begin{equation}\xlabel{ME}
\ME=\min\{\meh m,\,\#E-1\}.
\end{equation}

A function $f$ belongs to the trace space $\WMPR\vert  _E$ if and only if the following quantity
 \begin{equation}\xlabel{S-WP}
\TLNW(f:E)=
\,\left(\,\smed_{k=0}^{\ME}\,\,\smed_{i=\ell_1}^{\ell_2-k}
\min\left\{1,x_{i+m}-x_{i}\right\}
\,\left|\Delta^kf[x_i, \ldots ,x_{i+k}]\right|^p
\right)^{\frac1p}
\end{equation}
 is finite. (Note that, according to the notational convention adopted in \rf{AGR}, it must be understood that in \rf{S-WP} we have $x_i=+\infty$ for every $i>\ell_2$.)

Furthermore, $f$ belongs to $\WMIR\vert  _E$ if and only if the quantity
\[
\TLNWI(f:E)=\,\sup_{{k=0, \ldots ,\ME}\atop
{i=\ell_1, \ldots ,\ell_2-k}}
\,\,\,
\left|\Delta^kf[x_i, \ldots ,x_{i+k}]\right|
\]
 is finite.

Moreover, for every $p\in(1,\infty]$ the following equivalence $\| f\| _{\WMPR\vert  _E}\sim \TLNW(f:E)$ holds with constants depending only on $m$.
\end{theorem}

We prove Theorem \reff{W-TFIN} in Section 5. For a special case
of Theorem \reff{W-TFIN} for $m=2$ and strictly increasing
sequences $\{x_i\}_{i\in\Z}$ with  $x_{i+1}-x_i\le const$, and also for other
results related to characterization of the trace space
$\WMPR\vert  _E$, see Est\'evez~\cite{Es}.
\smallskip

In Section 6 we also give another characterization of the trace space $\WMPR\vert  _E$ expressed in terms of $L^p$-norms of certain kinds of ``sharp maximal functions'' which are defined as follows.
\begin{definition}
For each $m\in\N$, each closed set $E\subset\R$, each function $f:E\to\R$ and each integer $k$ in the range $0\le k \le m-1$, we let $\fks$ denote the maximal function associated with $f$ which is given by
 \begin{equation}\xlabel{FK-1}
\fks(x)=\sup_{{S\subset E,\,\#S=k+1,}\atop
{\diam(S\cup\{x\})\,\le 1}} \,\vert  \Delta^kf[S]\vert  ,~~~~~~x\in\R,
\end{equation}
 and, when $k=m$, by
 \begin{equation}\xlabel{FK-2}
\fms(x)=\sup_{{S\subset E,\,\#S=k+1, \vspace*{0.25mm}\atop
\diam(S\cup\{x\})\,\le 1}}
\,\frac{\diam S}{\diam(S\cup\{x\})} \,\vert  \Delta^mf[S]\vert  ,~~~~~~x\in\R.
\end{equation}

Here for every $k=0, \ldots ,m$, the above two suprema are taken over all $(k+1)$-point subsets $S\subset E$ such that $\diam(S\cup\{x\})\le 1$. If the family of sets $S$ satisfying these conditions is empty, we put $\fks(x)=0$.
\end{definition}

\begin{theorem}\lbl{W-MF}
Let $m\in\N$ and let $p\in(1,\infty)$. A function $f$ belongs to the trace space $\WMPR\vert  _E$ if and only if the function $\fks$ belongs to $\LPR$~for every $k=0, \ldots ,m$. Furthermore,
\[
\| f\| _{\WMPR\vert  _E}\sim \smed_{k=0}^m\,\,\| \fks\| _{\LPR}
\]
 with constants in this equivalence depending only on $m$ and $p$.
\end{theorem}

Here, as in our earlier paper~\cite{Sh-LMP-2018} we again feel a
strong debt to the remarkable papers of Calder\'{o}n and
Scott~\cite{C1,CS} which are devoted to characterization of
Sobolev spaces on $\RN$ in terms of  classical sharp
maximal functions. These papers motivated us to formulate and
subsequently prove Theorem \reff{W-MF} (as well as similar
results in~\cite{Sh-LMP-2018}). For an analog of
Theorem \reff{W-MF} for the space $\WORN$, $p>n$, we refer
the reader to~\cite[Theorem 1.4]{Sh2}.
     \msk

The next theorem states that there exists a solution to Problem \reff{PR-MAIN} which depends \textit{li\-nearly} on the initial data, i.e.,~the functions defined on $E$.
\begin{theorem}
 \lbl{LIN-OP} For every closed subset $E\subset\R$, every $p\in(1,\infty]$ and every $m\in\N$ there exists a continuous linear extension operator which maps the trace space $\WMPR\vert  _E$ into $\WMPR$. Its operator norm is bounded by a constant depending only on $m$.
\end{theorem}

As we have noted above, a variant of Problem \reff{PR-MAIN} for
the homogeneous Sobolev space $\LMPR$ has been studied in
the author's paper~\cite{Sh-LMP-2018}. In that paper we give
intrinsic characterizations of the restriction of
$\LMPR$-functions to an arbitrary closed subset
$E\subset\R$, and prove a series of $L^m_p$-extension
criteria similar to these given in Theorems 
\reff{W-VAR-IN}, \reff{W-TFIN} and \reff{W-MF}.

In particular, the variational extension criterion for $L^m_p$-traces given in Theorem \reff{MAIN-TH}, is one of the main ingredients of the proof of the sufficiency part of Theorem \reff{W-VAR-IN}. (Correspondingly, the proof of Theorem \reff{W-TFIN} relies on Theorem \reff{DEBOOR}, a counterpart of Theorem \reff{W-TFIN} for the space $\LMPR$.)
\smsk

Let us recall something of the history of extension problems for
the space $\LMPR$, and give a sketch of the proof of
Theorems \reff{W-VAR-IN}, \reff{W-TFIN} and \reff{LIN-OP}.
(See~\cite{Sh-LV-2018,Sh-LMP-2018} and references therein for
more details.)

Whitney~\cite{W2} completely solved an analog of
Problem \reff{PR-MAIN} for the space $C^m(\R)$. Whitney's extension
construction~\cite{W2} produces a certain \textit{extension
operator}
 \begin{equation}\xlabel{EXT-L}
\WOP:C^m(\R)\vert  _E\to C^m(\R)
\end{equation}
 which linearly and continuously maps the trace space  $C^m(\R)\vert  _E$ into $C^m(\R)$.

Whitney's extension method~\cite{W2} also provides a complete
solution to Problem \reff{PR-MAIN} for the space $\LMIR$.
Recall that $\LMIR$ can be identified with the space
$C^{m-1,1}(\R)$ of all $C^{m-1}$-functions on $\R$ whose
derivatives of order $m-1$ satisfy a Lipschitz condition.
In particular, the method of proof and technique developed
in~\cite{W2} (see also~\cite{Mer}) lead us to the following well
known description of the trace space $\LMIR\vert  _E$: A function
$f\in\LMIR\vert  _E$ if and only if the following quantity
 \begin{equation}\xlabel{L-INF1}
\NIN_{m,\infty}(f:E)=
\sup_{S\subset E,\,\,\# S=m+1}
\vert  \Delta^mf[S]\vert
\end{equation}
 is finite. Furthermore,
 \begin{equation}\xlabel{T-R1}
C_1\, \NIN_{m,\infty}(f:E)\le \| f\| _{\LMIR\vert  _E}\le C_2\, \NIN_{m,\infty}(f:E)
\end{equation}
 where $C_1$ and $C_2$ are positive constants depending only on $m$.

Favard~\cite{Fav} ($p=\infty$) and de Boor~\cite{deB4}
($p\in (1,\infty)$) characterized the traces of $\LMPR$-functions
to arbitrary sequences of points in $\R$. (Note that the
extension methods used in~\cite{Fav} and~\cite{deB4} are very
different from the Whitney extension construction~\cite{W2}.)
\begin{theorem}
\lbl{DEBOOR} Let $p\in(1,\infty]$, and let $\ell_1,\ell_2\in\Z\cup\{\pm\infty\}$, $\ell_1+m\le\ell_2$. Let $f$ be a function defined on a strictly increasing sequence of points $E=\{x_i\}_{i=\ell_1}^{\ell_2}$. We let $\TNMP(f:E)$ denote the quantity
\[
\TNMP(f:E)=\,\left(\smed_{i=\ell_1}^{\ell_2-m}\,\,
(x_{i+m}-x_i)\,\vert  \Delta^mf[x_i, \ldots ,x_{i+m}]\vert  ^p
\right)^{\frac1p}
\]
 for each $p\in(1,\infty)$, and
\[
\TNMI(f:E)=\sup_{\ell_1\le i\le  \ell_2-m}
\vert  \Delta^mf[x_i, \ldots ,x_{i+m}]\vert  .
\]

In these settings, the function $f$ is an element of $\LMPR\vert  _E$ if and only if $\TNMP(f:E)<\infty$.

Furthermore, $\| f\| _{\LMPR\vert  _E}\sim \TNMP(f:E)$ and the constants in this equivalence depend only on $m$.
\end{theorem}

Let us mention that in~\cite{Sh-LMP-2018} we have shown that
\textit{the very same Whitney extension operator $\WOP$}
(which was introduced in~\cite{W2} for characterization of the
trace space $C^{m}(\R)\vert  _E$) is ``universal'' for the scale of
$\LMPR$ spaces, $1<p\le\infty$. This means the following:
\textit{for every $p\in(1,\infty]$ the extension operator
$\WOP$ provides almost optimal $L^m_p$-extensions of
functions defined on $E$}. Furthermore, the operator
norm of $\WOP$ is bounded by a constant depending only on
$m$. That result  is an essential ingredient in our
proof in~\cite{Sh-LMP-2018} of the following generalization of
Theorem \reff{DEBOOR}:
\begin{theorem}
\lbl{MAIN-TH}(Variational Extension Criterion for
$\LMPR$-traces,~\cite{Sh-LMP-2018}) Let $p\in(1,\infty)$ and let
$m$ be a positive integer. Let $E\subset\R$ be a closed
set containing at least $m+1$ points. A function
$f:E\to\R$ can be extended to a function $F\in\LMPR$ if and
only if the following quantity
 \begin{equation}\xlabel{NMP}
\NMP(f:E)=\,\sup_{\{x_0, \ldots ,x_n\}\subset E}
\,\,\,\left(\smed_{i=0}^{n-m}
\,\,
(x_{i+m}-x_i)\,\vert  \Delta^mf[x_i, \ldots ,x_{i+m}]\vert  ^p
\right)^{\frac1p}
\end{equation}
 is finite. Here the supremum is taken over all integers $n\ge m$ and all strictly increasing sequences $\{x_0, \ldots ,x_n\}\subset E$ of $n+1$ elements. Furthermore,
\[
\| f\| _{\LMPR\vert  _E}\sim \NMP(f:E)
\]
 where the constants in this equivalence depend only on $m$.
\end{theorem}

For a more detailed account of the history of extension and
restriction problems for the Sobolev spaces $\LMPR$ and
$\WMPR$, and various references related to this topic, we
refer the reader to~\cite{Sh-LV-2018,Sh-LMP-2018}.

Let us now briefly describe the main ideas of the proof of the sufficiency in Theorem \reff{W-VAR-IN}, the most technically difficult part of the present paper. Let
\[
G=\{x\in\R: \dist(x,E)\ge 1\}~~~\text{and}~~~
\tE=E\cup G.
\]

Given a function $f$ on $E$, we let $\tf:\tE\to\R$ denote the extension of $f$ from $E$ to $\tE$ \textit{by zero}. (Thus, $\tf\vert  _{E}=f$ and $\tf\equiv 0$ on $\tE\setminus E$.) We show that $\tf$ satisfies the hypothesis of Theorem \reff{MAIN-TH} provided $f$ satisfies the hypothesis of Theorem \reff{W-VAR-IN} (i.e.,~$\NWMP(f:E)<\infty$). Our proof of this fact relies on a series of auxiliary lemmas devoted to $p$-summations of divided differences and connections between divided differences of various orders. See Section 3.2 and Section 4.

Next, Theorem \reff{MAIN-TH} tells us that $\tf$ can be extended to a function $F=\WOPTE(\tf)\in\LMPR$. Moreover, we prove that $\| F\| _{\LMPR}\le C_1(m)\,\NWMP(f:E)$. We also show that
\[
F\in\LPR~~~\text{and}~~~\| F\| _{\LPR}\le C_2(m)\,\NWMP(f:E).
\]

Hence we conclude that $F\in\WMPR$ and $\| F\| _{\WMPR}\le C_3(m)\,\NWMP(f:E)$, completing the proof of the sufficiency in Theorem \reff{W-VAR-IN}. (Here $C_i(m)$, $i=1,2,3$, are positive constants depending only on $m$.)
\smsk

In the same fashion we prove the sufficiency part of Theorem \reff{W-TFIN} replacing in this scheme Theorem \reff{MAIN-TH} with its refinement for sequences given in Theorem \reff{DEBOOR}.
\smsk

Furthermore, because $\tf$ depends \textit{linearly} on $f$ and the Whitney extension operator $\WOPTE$ is \textit{linear}, the extension operator $f\to F=\WOPTE(\tf)$ depends \textit{linearly} on $f$ as well. The necessity part of Theorem \reff{W-VAR-IN} (see Section 3.1) tells us that $\NWMP(f:E)\le C_4(m)\,\| f\| _{\WMPR\vert  _E}$. Hence,
\[
\| F\| _{\WMPR}\le C_3(m)\,\NWMP(f:E)
\le C_3(m)\cdot C_4(m)\,\| f\| _{\WMPR\vert  _E}.
\]

Thus, the mapping $f\to \WOPTE(\tf)$ provides a linear extension operator from $\WMPR\vert  _E$ into $\WMPR$ whose  operator norm is bounded by $C_3(m)\cdot C_4(m)$. This completes the proof of Theorem \reff{LIN-OP}.
\msk

Readers might find it helpful to also
consult the paper~\cite{Sh-LV-2018} posted on the arXiv.
It contains a much more detailed version of the
material presented here and also of the material presented
in~\cite{Sh-LMP-2018}.

\SECT{2. Divided differences: main properties.}{2}
\addtocontents{toc}{2. Divided differences: main properties. \hfill \thepage\par\VST}

\indent\par Let us fix some notation. Throughout the paper $C,C_1,C_2, \ldots $ will be generic posi\-tive constants which depend only on $m$ and $p$. These symbols may denote different constants in different occurrences. The dependence of a constant on certain parameters is expressed by the notation $C=C(m)$, $C=C(p)$ or $C=C(m,p)$. Given constants $\alpha,\beta\ge 0$, we write $\alpha\sim \beta$ if there is a constant $C\ge 1$ such that $\alpha/C\le \beta\le C\,\alpha$.

Given a measurable set $A\subset \R$, we let  $\left|A\right|$ denote the Lebesgue measure of $A$.
If $A\subset\R$ is finite, by $\#A$ we denote the number of elements of $A$. Given $A,B\subset \R$, let
\[
\diam A=\sup\{\mid a-a'\,\vert  :~a,a'\in A\}~~~\text{and}~~~\dist(A,B)=\inf\{\mid a-b\,\vert  :~a\in A, b\in B\}.
\]
 For $x\in \R$ we also set $\dist(x,A)=\dist(\{x\},A)$. Finally, we put $\dist(A,\emp)=+\infty$ provided  $A\ne\emp$.

Given $M>0$ and a family $\Ic$  of intervals in $\R$ we say that \textit{covering multiplicity} of $\Ic$ is bounded by $M$ if every point $x\in\R$ is covered by at most $M$ intervals from $\Ic$.

Given a function $g\in L_{1,loc}(\R)$ we let $\Mc[g]$ denote the Hardy--Littlewood maximal function of $g$:
 \begin{equation}\xlabel{HL-M}
\Mc[g](x)=\sup_{I\ni x}\frac{1}{\vert  I\vert  }\intl_I\vert  g(y)\vert  dy,~~~~x\in\R.
\end{equation}
 Here the supremum is taken over all closed intervals $I$ in $\R$ containing $x$.

By $\PM$ we denote the space of all polynomials of one variable of degree at most $m$ defined on $\R$. Finally, given a nonnegative integer $k$, a $(k+1)$-point set $S\subset\R$ and a function $f$ on $S$, we let $L_S[f]$ denote the Lagrange polynomial of degree at most $k$ interpolating $f$ on $S$; thus
$L_S[f]\in \Pc_k$ and $L_S[f](x)=f(x)$ for every $x\in S$.
\smallskip

In this section we recall a definition and several useful
properties of the divided differences of functions. We refer the
reader to~\cite[Ch. 4, \textsection 7]{DL},~\cite[Section 1.3]{FK-88}
and~\cite[Section 2.1]{Sh-LV-2018} for the proofs of these
properties.
 \msk

First, we recall the standard definition of the divided difference. Let $k$ be a nonnegative integer, and let $S=\{x_0, \ldots ,x_k\}$ be a $(k+1)$-point subset of $\R$.

We set $\Delta^0f[S]=f(x_0)$ provided $S=\{x_0\}$ is a singleton, and given $k\in\N$ we set
 \begin{align}\xlabel{D-IND}
\Delta^kf[S]&=\Delta^{k}f[x_0,x_1, \ldots ,x_{k}]\nonumber\\
&=\left(\Delta^{k-1}f[x_1, \ldots ,x_{k}]
-\Delta^{k-1}f[x_0, \ldots ,x_{k-1}]\right)/(x_k-x_0).
\end{align}

This definition implies the following very well known formula for the divided difference:
 \begin{equation}\xlabel{D-PT1}
\Delta^kf[S]=
\smed_{i=0}^k\,\,\frac{f(x_i)}{\omega'(x_i)}=
\smed_{i=0}^k\,\,\frac{f(x_i)}
{\prod_{j\in\{0, \ldots ,k\},j\ne i}(x_i-x_j)}
\end{equation}
 where $\omega(x)=(x-x_0)...(x-x_k)$ and $\omega'$ is the derivative of $\omega$.

We recall that $L_S[f]$ denotes the Lagrange polynomial of degree at most $k=\#S-1$ interpolating $f$ on $S$. Then the following equality
 \begin{equation}\xlabel{D-LAG}
\Delta^{k}f[S]=\frac{1}{k!}\,L^{(k)}_S[f]
\end{equation}
 holds. Thus, $\Delta^{k}f[S]=A_k$ where $A_k$ is the coefficient of $x^k$ in the polynomial $L_S[f]$.
\smallskip

In Sections 3, 4 and 6 we will need the following properties of divided differences:
\smsk

$(\bigstar 1)$~Let $k\in\N$, and let $x_0=\min\{x_i:i=0, \ldots ,k\}$ and $x_k=\max\{x_i:i=0, \ldots ,k\}$. Then for every function $F\in C^k[x_0,x_k]$ there exists $\xi\in [x_0,x_k]$ such that
 \begin{equation}\xlabel{D-KSI}
\Delta^{k}F[x_0,x_1, \ldots ,x_{k}]
=\frac{1}{k!}\,F^{(k)}(\xi)\,.
\end{equation}

$(\bigstar 2)$~Let $x_0<x_1<\cdots<x_m$, and let $F$ be a function on $[x_0,x_m]$ with absolutely continuous derivative of order $m-1$. Then
 \begin{equation}\xlabel{DVD-IN}
\vert  \Delta^{m}F[S]\vert  \le\frac{1}{(m-1)!}\cdot\frac{1}{x_m-x_0}
\,\intl_{x_0}^{x_m}\,
\vert  F^{(m)}(t)\vert  \,dt.
\end{equation}

$(\bigstar 3)$~(\cite[p. 15]{FK-88}) Let $k,n\in\N$,
$n\ge k$, and let $\{s_0, \ldots ,s_{n}\}\subset\R$, be a strictly increasing
sequence. Let $g$ be a function on $\{s_0, \ldots ,s_{n}\}$, and let
$\{t_0, \ldots ,t_k\}$ be a $(k+1)$-point subset of $\{s_0, \ldots ,s_{n}\}$.

There exist $\alpha_i\in\R$, $\alpha_i\ge 0$, $i=0, \ldots ,n-k$, such that $\alpha_0+\cdots+\alpha_{n-k}=1$ and
\[
\Delta^kg[t_0, \ldots ,t_k]=\smed_{i=0}^{n-k}\,\,\alpha_i\,
\Delta^kg[s_i, \ldots ,s_{i+k}].
\]


\SECT{3. The variational extension criterion: auxiliary results.}{3}
\addtocontents{toc}{3. The variational extension criterion: auxiliary results. \hfill\thepage\par \VST}

\indent\par In this section we prove a series of auxiliary lemmas which we will need for the proofs of Theorem \reff{W-VAR-IN} (see Sections 4) and Theorem \reff{W-TFIN} (see Sections 5).
\bsk

\par {\bf 3.1. The necessity part of the variational extension criterion.}
\medskip
\addtocontents{toc}{~~~~3.1. The necessity part of the variational extension criterion. \hfill \thepage\par}

Let $m$ be a positive integer and let $E$ be a closed subset of $\R$ containing at least $m+1$ points.

We begin the proof of the necessity part of Theorem \reff{W-VAR-IN} with the case $p\in(1,\infty)$.

Let $f\in\WMPR\vert  _E$, and let $F\in\WMPR$ be a function such that $F\vert_E=f$. Let $n\ge m$ and let $\{x_0, \ldots ,x_n\}$ be a finite strictly increasing sequences in $E$. Our aim is to prove that for each $k=0, \ldots ,m$, the following inequality
 \begin{equation}\xlabel{S-W1}
L_k=\smed_{i=0}^{n-k}
\min\left\{1,x_{i+m}-x_{i}\right\}
\,\left|\Delta^kf[x_i, \ldots ,x_{i+k}]\right|^p
\le C(m)^p\,\| F\| ^p_{\WMPR}
\end{equation}
 holds.

First, prove this inequality for $k=m$.  Thanks to \rf{DVD-IN} and H\"{o}lder's inequality, for every $i=0, \ldots ,n-m$, we have
\begin{eqnarray}
A_i&=&(x_{i+m}-x_i)\,\vert  \Delta^m f[x_i, \ldots ,x_{i+m}]\vert  ^p
=
(x_{i+m}-x_i)\,\vert  \Delta^mF[x_i, \ldots ,x_{i+m}]\vert  ^p\nn\\
&\le&
(x_{i+m}-x_i)\cdot
\left\{\frac{1}{(m-1)!}\cdot
\,\frac{1}{x_{i+m}-x_i}
\,\intl_{x_i}^{x_{i+m}}\,\vert  F^{(m)}(t)\vert  \,dt\right\}^p\nonumber\\
&\le&\,
\frac{1}{((m-1)!)^p}
\,\,\intl_{x_i}^{x_{i+m}}\,\vert  F^{(m)}(t)\vert  ^p\,dt\,.\nn
\end{eqnarray}
 Hence,
\[
L_m\le \smed_{i=0}^{n-m}
\,\,A_i
\,\le
\frac{1}{((m-1)!)^p}\,\smed_{i=0}^{n-m}
\,\,\intl_{x_i}^{x_{i+m}}\,\vert  F^{(m)}(t)\vert  ^p\,dt.
\]

Clearly, the covering multiplicity of the family $\{(x_i,x_{i+m}):i=0, \ldots ,n-m\}$ of open intervals is bounded by $m$, so that
\[
L_m\le \,\frac{m}{((m-1)!)^p}
\,\intl_{x_0}^{x_{n}}\,\vert  F^{(m)}(t)\vert  ^p\,dt\le 2^p\,\| F\| _{\LMPR}^p\le 2^p\,\| F\| ^p_{\WMPR}
\]
 proving \rf{S-W1} for $k=m$.

Let us prove \rf{S-W1} for  $k\in\{0, \ldots ,m-1\}$. To do so, we need the following
\begin{lemma}
\lbl{NP-W} Let  $q\in[1,\infty)$. Let $k\in\{0, \ldots ,m-1\}$, and let $S$ be a $(k+1)$-point subset of a closed interval $I\subset\R$ (bounded or unbounded). Let $G\in C^k[I]$ and let $G^{(k)}$ be absolutely continuous on $I$ with $G^{(k+1)}\in L_q[I]$. Then the following inequality
 \begin{equation}\xlabel{T-F1}
\min\{1,\vert  I\vert  \}\cdot\vert  \Delta^{k}G[S]\vert  ^q\le 2^q\,\intl_I \left(\vert  G^{(k)}(y)\vert  ^q+\vert  G^{(k+1)}(y)\vert  ^q\right)\,dy
\end{equation}
 holds.
\end{lemma}

{\it Proof.} Let $S=\{y_0, \ldots ,y_k\}$ with $y_0<\cdots<y_k$. Property \rf{D-KSI} tells us that there exists $s\in[y_0,y_k]$ such that
$k!\,\Delta^{k}G[S]=G^{(k)}(s).$

Let $t=\min\{1,\vert  I\vert  \}$. Because $s\in I$ and $t\le \vert  I\vert  $, there exists a closed interval $J\subset I$ with $\vert  J\vert  =t$ such that $s\in J$. It now follows that for every $y\in J$
\begin{align*}
\vert  \Delta^{k}G[S]\vert  ^q&\le \vert  G^{(k)}(s)\vert  ^q\le 2^q\,\vert  G^{(k)}(s)-G^{(k)}(y)\vert  ^q+2^q\,\vert  G^{(k)}(y)\vert  ^q\nonumber\\
&\le 2^q\left(\intl_J \vert  G^{(k+1)}(x)\vert  \,dx\right)^q+2^q\,\vert  G^{(k)}(y)\vert  ^q.
\end{align*}

This and H\"{o}lder's inequality imply that
\[
\vert  \Delta^{k}G[S]\vert  ^q\le 2^q\,\vert  J\vert  ^{q-1}\,
\intl_J \vert  G^{(k+1)}(x)\vert  ^q dx+2^q\,\vert  G^{(k)}(y)\vert  ^q,~~~y\in J.
\]
 Integrating this inequality on $J$ with respect to $y$, we obtain the following:
\[
\vert  J\vert  \,\vert  \Delta^{k}G[S]\vert  ^q\le 2^q\,\vert  J\vert  ^{q}
\intl_J \vert  G^{(k+1)}(x)\vert  ^q dx+2^q\intl_J\vert  G^{(k)}(y)\vert  ^q\,dy.
\]

Because $\vert  J\vert  =t=\min\{1,\vert  I\vert  \}\le 1$ and $J\subset I$, this inequality implies \rf{T-F1} proving the lemma.\bx

Fix $k\in\{0, \ldots ,m-1\}$. Let $i\in\{0, \ldots ,n-k\}$,  $T_i=[x_{i},x_{i+m}]$ and
\[
t_i=\min\{1,\vert  T_i\vert  \}=\min\{1,x_{i+m}-x_{i}\}.
\]
 (We recall the convention adopted in \rf{AGR} which tells us that $x_{i+m}=+\infty$ if $i+m>n$.)

Let
\[
B_i=
\min\left\{1,x_{i+m}-x_{i}\right\}
\,\left\lvert  \Delta^kf[x_i, \ldots ,x_{i+k}]\right\rvert  ^p.
\]

Let us apply Lemma \reff{NP-W} to $G=F$, $t=t_i$, $q=p$, and $S=\{x_i, \ldots ,x_{i+k}\}$. (Clearly, $F^{(k+1)}\in L_p(\R)$ because $F\in\WMPR$ and $k\le m-1$). This lemma tells us  that
\[
B_i\le 2^p\,\intl_{T_i}
\left(\vert  F^{(k)}(y)\vert  ^p+\vert  F^{(k+1)}(y)\vert  ^p\right)\,dy.
\]
 Hence,
\[
L_k=\smed_{i=0}^{n-k}\,B_i\le 2^p\,\smed_{i=0}^{n-k}\,\intl_{T_i}
\left(\vert  F^{(k)}(y)\vert  ^p+\vert  F^{(k+1)}(y)\vert  ^p\right)\,dy.
\]

Note that the covering multiplicity of the family $\{T_i=[x_i,x_{i+m}]:i=0, \ldots ,n-k\}$ of closed intervals is bounded by $m+1$. Hence,
\[
L_k\le
2^p\,(m+1)\intl_{\R}\,
\left(\vert  F^{(k)}(y)\vert  ^p+\vert  F^{(k+1)}(y)\vert  ^p\right)\,dy
\le C(m)^p\,\| F\| ^p_{\WMPR}
\]
 proving \rf{S-W1} for every integer $k$ in the range $0\le k\le m-1$.
\smsk

Thus \rf{S-W1} holds for all $k=0, \ldots ,m$. This inequality and definition \rf{N-WRS-IN} imply that
 \begin{equation}\xlabel{NW-F1}
\NWMP(f:E)\le C(m)\,\| F\| _{\WMPR}
\end{equation}
 for every $F\in\WMPR$ such that $F\vert_E=f$. Hence,
 \begin{equation}\xlabel{NW-F2}
\NWMP(f:E)\le C(m)\,\| f\| _{\WMPR\vert  _E}~~~~~\text{(see \rf{N-WMPR}),}
\end{equation}
   completing the proof of the necessity part of Theorem \reff{W-VAR-IN} for the case $1<p<\infty$.
\smallskip

Prove the necessity for $p=\infty$ and arbitrary closed set $E$ (not necessarily with $\#E\ge m+1$).

Property $(\bigstar 2)$ of Sections 2 tells us that if $F\in\WMIR$ and $F\vert_E=f$, then for every $k\in\{0, \ldots ,m\}$ and any strictly increasing sequence $\{y_0, \ldots ,y_k\}\subset E$ the following inequality
\[
\vert  \Delta^{k}f[y_0, \ldots ,y_k]\vert  =
\vert  \Delta^{k}F[y_0, \ldots ,y_k]\vert  \le\frac{1}{(k-1)!}
\cdot\frac{1}{y_k-y_0}
\,\intl_{y_0}^{y_k}\,\vert  F^{(k)}(t)\vert  \,dt
\]
 holds. Hence, $\vert  \Delta^{k}f[y_0, \ldots ,y_k]\vert  \le \| F\| _{\WMIR}$ proving that $\NWMI(f:E)\le \| F\| _{\WMIR}$, see \rf{N-WINF}.

Thus, \rf{NW-F1} holds for $p=\infty$ proving  \rf{NW-F2} for this case.

The proof of the necessity part of Theorem \reff{W-VAR-IN} is complete.\bx
\bsk

\par {\bf 3.2. Preparatory results for the sufficiency part of the variational extension criterion.}
\medskip
\addtocontents{toc}{~~~~3.2. Preparatory results for the sufficiency part of the variational extension criterion. \hfill \thepage\par\VST}

In this section we adopt the following extended version of the convention adopted in \rf{AGR}:
\begin{convention}
\lbl{AGREE-2} Given a finite strictly increasing sequence $\{y_i\}_{i=0}^n\subset \R$ we put
\[
y_j=-\infty~~~~\text{if}~~~j<0,~~~~\text{and}~~~~
y_j=+\infty~~~\text{if}~~~j>n~~~~\text{(as in \rf{AGR})}.
\]
\end{convention}

Our proof of the sufficiency part of Theorem \reff{W-VAR-IN} relies on a series of auxiliary lemmas.
The first of them is the following well known statement from graph theory.
\begin{lemma}
\lbl{GRAPH} Let $\ell\in\N$ and let $\Ac=\{A_\alpha:\alpha\in I\}$ be a family of subsets of\, $\R$ such that every set $A_\alpha\in\Ac$ has common points with at most $\ell$ sets $A_\beta\in\Ac$.

Then there exist subfamilies $\Ac_i\subset \Ac$, $i=1, \ldots ,n$, with $n\le \ell+1$, each consisting of pairwise disjoint sets such that
\[
\Ac=\bigcup_{i=1}^n\,\Ac_j.
\]
\end{lemma}

{\it Proof.} The proof is immediate from the following well-known
statement (see, e.g.~\cite{JT}): \textit{Every graph can be
colored with one more color than the maximum vertex degree.}\bx

\begin{lemma}
\lbl{SM-DF} Let $n\in\N$ and let $S=\{y_0,\ldots,y_n\}$ where $\{y_i\}_{i=0}^n$ is a strictly increasing sequence of points in $\R$ such that $\diam S=y_n-y_0\ge 1$. Let $h$ be a function on $S$.

Then there exist $k\in \{0, \ldots ,n-1\}$ and $i\in\{0, \ldots ,n-k\}$ such that $y_{i+k}-y_i\le 1$ and
 \begin{equation}\xlabel{Y-L}
\vert  \Delta^nh[S]\vert  \le\,2^n \,\vert  \Delta^kh[y_i, \ldots ,y_{i+k}]\vert  /\diam S\,.
\end{equation}
 Furthermore,
 \begin{equation}\xlabel{SV-T}
\text{either}~~~i+k+1\le n~~\text{and}~~
y_{i+k+1}-y_{i}\ge 1,~~\text{or}~~i\ge 1 ~~\text{and}~~
y_{i+k}-y_{i-1}\ge 1\,.
\end{equation}
\end{lemma}

{\it Proof.} We proceed by induction on $n$. Let $n=1$, and let $S=\{y_0,y_1\}$ where $y_0<y_1$ and $y_1-y_0\ge 1$. Then for every function $h$ on $S$ the following inequality
\[
\vert  \Delta^{1}h[y_0,y_1]\vert  =\frac{\vert  h(y_1)-h(y_0)\vert  }{y_1-y_0}\le \frac{\vert  h(y_1)\vert  +\vert  h(y_0)\vert  }{y_1-y_0}\le \frac{2\max\{\vert  h(y_0)\vert  ,\vert  h(y_1)\vert  \}}{y_1-y_0}
\]
 holds. Let us pick $i\in\{0,1\}$ such that $\vert  h(y_i)\vert  =\max\{\vert  h(y_0)\vert  ,\vert  h(y_1)\vert  \}$. Then
\[
\vert  \Delta^{1}h[y_0,y_1]\vert  \le
\frac{2\max\{\vert  \Delta^{0}h[y_0]\vert  ,\vert  \Delta^{0}h[y_1]\vert  \}}
{y_1-y_0}=
2\,\vert  \Delta^{0}h[y_i]\vert  /(y_1-y_0)
\]
 proving \rf{Y-L} for $n=1$ with $k=0$. It is also clear that $i+k+1\le n$ and $y_{i+k+1}-y_{i}\ge 1$ if $i=0$, and $i\ge 1$ and $y_{i+k}-y_{i-1}\ge 1$ if $i=1$ proving  \rf{SV-T} and the lemma for $n=1$.

For the induction step, we fix $n\ge 1$ and suppose the lemma holds for $n$; we then prove it for $n+1$.

Let $\{x_i\}_{i=0}^{n+1}$ be a strictly increasing sequence in $\R$ such that $x_{n+1}-x_0\ge 1$, and let $h:S\to\R$ be a function on the set $S=\{x_0,\ldots,x_{n+1}\}$. Then, thanks to \rf{D-IND},
\[
\Delta^{n+1}h[S]=\Delta^{n+1}h[x_0, \ldots ,x_{n+1}]
=\left(\Delta^{n}h[x_1, \ldots ,x_{n+1}]
-\Delta^{n}h[x_0, \ldots ,x_{n}]\right)/(x_{n+1}-x_0)
\]
 so that
 \begin{equation}\xlabel{T-33}
\vert  \Delta^{n+1}h[S]\vert
\le(\vert  \Delta^{n}h[x_1, \ldots ,x_{n+1}]\vert  +
\vert  \Delta^{n}h[x_0, \ldots ,x_{n}]\vert  )/\diam S.
\end{equation}

Let
\[
\SH_1=\{x_0, \ldots ,x_{n}\}~~~\text{and}~~~
\SH_2=\{x_1, \ldots ,x_{n+1}\}.
\]

If $\diam  \SH_1=x_n-x_0\ge 1$, then, thanks to the induction hypothesis, there exist $k_1\in\{0, \ldots ,n-1\}$ and $i_1\in\{0, \ldots ,n-k_1\}$ such that $x_{i_1+k_1}-x_{i_1}\le 1$ and
 \begin{equation}\xlabel{E-0}
\vert  \Delta^nh[\SH_1]\vert  \le\,2^n \,\vert  \Delta^{k_1}h[x_{i_1}, \ldots ,x_{i_1+k_1}]\vert  /\diam \SH_1\le
\,2^n \,\vert  \Delta^{k_1}h[x_{i_1}, \ldots ,x_{i_1+k_1}]\vert  \,.
\end{equation}
 Moreover,
 \begin{equation}\xlabel{E-1}
\text{either}~~i_1+k_1+1\le n~\text{and}~
x_{i_1+k_1+1}-x_{i_1}\ge 1,~\text{or}~~i_1\ge 1 ~\text{and}~
x_{i_1+k_1}-x_{i_1-1}\ge 1.
\end{equation}

Clearly, if $\diam  \SH_1\le 1$, then the inequality
$x_{i_1+k_1}-x_{i_1}\le 1$ and \rf{E-0} hold provided $k_1=n$ and $i_1=0$. It is also clear that in this case
\[
i_1+k_1+1\le n+1~~~\text{and}~~~~
x_{i_1+k_1+1}-x_{i_1}=x_{n+1}-x_0\ge 1.
\]

This observation and \rf{E-1} imply the following statement: If $i=i_1$ and $k=k_1$, then $k\in \{0, \ldots ,n\}$ and $i\in\{0, \ldots ,n+1-k\}$, and
\begin{equation}\xlabel{F-E3}
\text{either}~~i+k+1\le n+1~\text{and}~
x_{i+k+1}-x_{i}\ge 1,~\text{or}~~i\ge 1~\text{and}~
x_{i+k}-x_{i-1}\ge 1.
\end{equation}

In the same way we prove the existence of $k_2\in \{0, \ldots ,n\}$ and $i_2\in\{0, \ldots ,n+1-k_2\}$ such that $x_{i_2+k_2}-x_{i_2}\le 1$,
 \begin{equation}\xlabel{E-2}
\vert  \Delta^nh[\SH_2]\vert  \le\,2^n \,\vert  \Delta^{k_2}h[x_{i_2}, \ldots ,x_{i_2+k_2}]\vert  ,
\end{equation}
 and \rf{F-E3} holds provided $i=i_2$ and $k=k_2$.

Let us pick $\ell\in\{1,2\}$ such that
\[
\vert  \Delta^{k_{\ell}}h[x_{i_{\ell}}, \ldots ,x_{i_{\ell}+k_{\ell}}]\vert
=\max\{\vert  \Delta^{k_1}h[x_{i_1}, \ldots ,x_{i_1+k_1}]\vert  ,
\vert  \Delta^{k_2}h[x_{i_2}, \ldots ,x_{i_2+k_2}]\vert  \}.
\]
 Inequalities \rf{T-33}, \rf{E-0} and \rf{E-2} tell us that
\begin{align*}
\vert  \Delta^{n+1}h[S]\vert
&\le 2^n \,(\vert  \Delta^{k_1}h[x_{i_1}, \ldots ,x_{i_1+k_1}]\vert  +
\vert  \Delta^{k_2}h[x_{i_2}, \ldots ,x_{i_2+k_2}]\vert  )/\diam S\\
&\le 2^{n+1} \frac{\vert  \Delta^{k_{\ell}}
h[x_{i_{\ell}}, \ldots ,x_{i_{\ell}+k_{\ell}}]\vert  }
{\diam S}.
\end{align*}

Thus,
\[
\vert  \Delta^{n+1}h[S]\vert  \le\,2^{n+1}
\,\vert  \Delta^k h[x_i, \ldots ,x_{i+k}]\vert  /\diam S
\]
 provided $i=i_{\ell}$ and $k=k_{\ell}$. Furthermore, $k\in\{0, \ldots ,n\}$, $i\in \{0, \ldots ,n+1-k\}$, and the statement \rf{F-E3} holds for these $i$ and $k$. This proves the lemma for $n+1$ completing the proof.\bx

The next two lemmas are variants of results proven in~\cite{Es}.
\begin{lemma}\lbl{K-I1}
Let $m\in\N$, $p\in[1,\infty)$ and $k\in\{0, \ldots ,m-1\}$. Let $I\subset\R$ be a bounded interval, and let $z_0, \ldots ,z_{m-1}$ be $m$ distinct points in $I$. Then for every function $F\in\LMPR$ and every $x\in I$ the following inequality
\begin{equation}\xlabel{X-F}
\vert  F^{(k)}(x)\vert  ^p\le C^p
\left\{\vert  I\vert  ^{(m-k)p-1}\intl_I\vert  F^{(m)}(s)\vert  ^pds
+\smed_{j=k}^{m-1}
\vert  I\vert  ^{(j-k)p}\vert  \Delta^jF[z_0, \ldots ,z_j]\vert  ^p\right\}
\end{equation}
 holds. Here $C$ is a constant depending only on $m$.
\end{lemma}

{\it Proof.} Because $F\vert_I\in C^j[I]$ for every $j\in\{k, \ldots ,m-1\}$, property \rf{D-KSI} implies the existence of a point $y_j\in I$ such that
 \begin{equation}\xlabel{SI-L}
\Delta^{j}F[z_0, \ldots ,z_{j}]=\frac{1}{j!}\,F^{(j)}(y_j).
\end{equation}

The Newton--Leibniz formula tells us that
\[
F^{(j)}(y)=F^{(j)}(y_j)+\intl_{y_j}^y\,F^{(j+1)}(s)\,ds
~~~\text{for every}~~~y\in I.
\]
 Hence,
 \begin{equation}\xlabel{S-D}
\vert  F^{(j)}(y)\vert  \le \vert  F^{(j)}(y_j)\vert  +\intl_I\,\vert  F^{(j+1)}(s)\vert  \,ds,~~~y\in I.
\end{equation}

Because $x\in I$, we have
\[
\vert  F^{(k)}(x)\vert  \le \vert  F^{(k)}(y_k)\vert  +\intl_I\,\vert  F^{(k+1)}(s)\vert  \,ds.
\]
 Furthermore, \rf{S-D} implies that
\[
\vert  F^{(k+1)}(s)\vert  \le \vert  F^{(k+1)}(y_{k+1})\vert  +\intl_I\,\vert  F^{(k+2)}(t)\vert  \,dt~~~
\text{for every}~~~s\in I,
\]
 provided $k+2\le m$. Hence,
\[
\vert  F^{(k)}(x)\vert  \le \vert  F^{(k)}(y_k)\vert  + \vert  I\vert  \,\vert  F^{(k+1)}(y_{k+1})\vert  +\vert  I\vert  \,\intl_I\,\vert  F^{(k+2)}(t)\vert  \,dt.
\]

Repeating this inequality $m-k-1$ times, we obtain the following inequality:
\[
\vert  F^{(k)}(x)\vert  \le \smed_{j=k}^{m-1}
\,\vert  I\vert  ^{j-k}\,\vert  F^{(j)}(y_j)\vert  + \vert  I\vert  ^{m-k-1}\,\intl_I\,\vert  F^{(m)}(s)\vert  \,ds.
\]
 From this inequality and H\"{o}lder's inequality, we have
\[
\vert  F^{(k)}(x)\vert  \le \smed_{j=k}^{m-1}
\,\vert  I\vert  ^{j-k}\,\vert  F^{(j)}(y_j)\vert  + \vert  I\vert  ^{m-k-1}\,\vert  I\vert  ^{1-1/p}
\left(\intl_I\,\vert  F^{(m)}(s)\vert  ^p\,ds\right)^{1/p}
\]
 so that
\[
\vert  F^{(k)}(x)\vert  ^p\le (m-k+1)^{p-1}
\left\{\smed_{j=k}^{m-1}
\,\vert  I\vert  ^{(j-k)p}\,\vert  F^{(j)}(y_j)\vert  ^p+ \vert  I\vert  ^{(m-k)p-1}
\,\intl_I\,\vert  F^{(m)}(s)\vert  ^p\,ds\right\}.
\]
 This inequality and \rf{SI-L} imply \rf{X-F} proving the lemma.\bx

\begin{lemma}
\lbl{DF-CD} Let $Z=\{z_0, \ldots ,z_m\}$ be an $(m+1)$-point subset of $\R$, and let $g$ be a function on $Z$. Then for every $k=0, \ldots ,m-1$, and every $S\subset Z$ with $\#S=k+1$ the following inequality
\[
\vert  \Delta^kg[S]\vert  \le C(m)\,\smed_{j=k}^{m}
\,(\diam Z)^{j-k}\cdot\vert  \Delta^jg[z_0, \ldots ,z_j]\vert
\]
 holds.
\end{lemma}

{\it Proof.} Let $I=[\min Z,\max Z]$. Then $I\supset Z$ and $\vert  I\vert  =\diam Z$. Let $L_Z[g]$ be the Lagrange polynomial of degree at most $m$ which agrees with $g$ on $Z$. Thanks to \rf{D-LAG} and \rf{D-KSI},
\[
\Delta^{m}g[Z]=\frac{1}{m!}\,(L_Z[g])^{(m)}~~~
\text{and}~~~
\Delta^{k}g[S]=\frac{1}{k!}\,(L_Z[g])^{(k)}(\xi)~~~~\text{for some}~~~\xi\in I.
\]
 We apply Lemma \reff{K-I1} to the function $F=L_Z[g]$, $x=\xi$ and $p=1$, and obtain the following:
 \begin{eqnarray}
\vert  \Delta^{k}g[S]\vert
&\le&
C(m)\,
\left\{\vphantom{\sum_{j=k}^{m-1}}
\vert  I\vert  ^{(m-k)-1}\,\intl_I\,\vert  (L_Z[g])^{(m)}(s)\vert  \,ds\right.\nn\\
&&+ \left. \smed_{j=k}^{m-1}
\,\vert  I\vert  ^{j-k}\,\,\vert  \Delta^j(L_Z[g])[z_0, \ldots ,z_j]\vert  \right\}\nn\\
&\le&
C(m)\,
\left\{\vert  I\vert  ^{m-k}\,\vert  \Delta^{m}g[Z]\vert
+\smed_{j=k}^{m-1}
\,\vert  I\vert  ^{j-k}\,\,\vert  \Delta^jg[z_0, \ldots ,z_j]\vert  \right\}.
\nn
\end{eqnarray}

The proof of the lemma is complete. \bx

Integrating both sides of inequality \rf{X-F} on $I$ (with respect to $x$), we obtain the next
\begin{lemma}
\lbl{K-I2} In the settings of Lemma \reff{K-I1}, the following inequality
$$
\intl_I\vert  F^{(k)}(x)\vert  ^p\,dx \le C(m)^p\,
\left\{
\sum_{j=k}^{m-1}
\vert  I\vert  ^{(m-k)p}\,\intl_I\,\vert  F^{(m)}(s)\vert  ^p\,ds+\smed_{j=k}^{m-1}
\,\vert  I\vert  ^{(j-k)p+1}\,\,\vert  \Delta^jF[z_0, \ldots ,z_j]\vert  ^p\right\}
$$
holds.
\end{lemma}

\begin{lemma}
\lbl{P-WQ} Let $p\in[1,\infty)$ and $m,n\in\N$, $m\le n$. Let $S=\{y_0, \ldots ,y_n\}$ where $\{y_i\}_{i=0}^n$ is a strictly increasing sequence of points in $\R$. Suppose that there exists $\ell\in\N$, $m\le \ell\le n$, such that
\[
y_\ell-y_0\le 2~~~\text{but}~~~y_{\ell+1}-y_0>2.
\]

Then for every function $g$ defined on $S$ and every $k\in\{0, \ldots ,m-1\}$ the following inequality
 \begin{equation}\xlabel{L-M1}
\vert  \Delta^kg[y_0, \ldots ,y_k]\vert  ^p\le  C(m)^p\,
\smed_{j=k}^{m}\,\,\smed_{i=0}^{\ell-j}
\min\left\{1,y_{i+m}-y_{i}\right\}
\,\left\lvert  \Delta^jg[y_i, \ldots ,y_{i+j}]\right\rvert  ^p
\end{equation}
 holds. (We recall that, according to our convention \rf{AGR}, $y_i=+\infty$ provided $i>n$.)
\end{lemma}

{\it Proof.} Let $V=\{y_0, \ldots ,y_\ell\}$, and let
 \begin{equation}\xlabel{A-D}
A=\smed_{j=k}^{m}\,\,\smed_{i=0}^{\ell-j}
\min\left\{1,y_{i+m}-y_{i}\right\}
\,\left\lvert  \Delta^jg[y_i, \ldots ,y_{i+j}]\right\rvert  ^p.
\end{equation}
 Theorem \reff{DEBOOR} tells us that there exists a function $G\in\LMPR$ such that $G\vert_V=g\vert  _V$ and
\[
\| G\| ^p_{\LMPR}\le C(m)^p\,
\smed_{i=0}^{\ell-m}\,\,
(y_{i+m}-y_i)\,\vert  \Delta^mg[y_i, \ldots ,y_{i+m}]\vert  ^p.
\]

Note that $y_{i+m}-y_i\le y_\ell-y_0\le 2$ for $0\le i\le\ell-m$, so that
$y_{i+m}-y_i\le 2\min\{1,y_{i+m}-y_i\}$. Hence,
 \begin{equation}\xlabel{G-A1}
\| G\| ^p_{\LMPR}\le C(m)^p\,
\smed_{i=0}^{\ell-m}\,\,
\min\{1,y_{i+m}-y_i\}\,\vert  \Delta^mg[y_i, \ldots ,y_{i+m}]\vert  ^p\le C(m)^p\,A.
\end{equation}
 Furthermore, thanks to \rf{D-KSI}, there exists $\xb\in[y_0,y_k]$ such that
 \begin{equation}\xlabel{X-B}
\Delta^{k}g[y_0, \ldots ,y_{k}]
=\frac{1}{k!}\,G^{(k)}(\xb)\,.
\end{equation}

Let $I=[y_0,y_\ell]$, and let $V_k=\{y_0, \ldots ,y_k\}$, $0\le k\le \ell$. Consider the following two cases.

\smallskip
\textit{The first case:} $y_\ell-y_0\le 1$.~Let
 \begin{equation}\xlabel{Z-N}
z_j=y_{\ell-m+1+j},~~~~j=0, \ldots ,m-1.
\end{equation}

Because $\xb, z_0, \ldots ,z_{m-1}\in I$, property \rf{X-B} and Lemma \reff{K-I1} imply that
\begin{equation}
\xlabel{DK-11}
\vert  \Delta^{k}g[V_k]\vert  ^p
\le C^p
\left\{
\vert  I\vert  ^{(m-k)p-1}\intl_I\vert  G^{(m)}\vert  ^pds +\smed_{j=k}^{m-1}
\vert  I\vert  ^{(j-k)p}\vert  \Delta^jG[z_0, \ldots ,z_j]\vert  ^p\right\}
\end{equation}
where $C$ is a constant depending only on $m$.

Note that $\vert  I\vert  =y_\ell-y_0\le 1$ and $(m-k)p-1\ge 0$ because $0\le k\le m-1$ and $p\ge 1$. Therefore, $\vert  I\vert  ^{(m-k)p-1}\le 1$ and $\vert  I\vert  ^{(j-k)p}\le 1$ for every $j=k, \ldots ,m-1$. From this and inequalities  \rf{G-A1}, \rf{Z-N} and \rf{DK-11}, it follows that
\[
\vert  \Delta^{k}g[V_k]\vert  ^p
\le C(m)^p\,
\left\{A+\smed_{j=k}^{m-1}
\,\vert  \Delta^jg[y_{\ell-m+1}, \ldots ,y_{\ell-m+1+j}]\vert  ^p\right\}.
\]

Because $y_\ell-y_0\le 1$ and $y_{\ell+1}-y_0>2$, we have $\min\{1,y_{\ell+1}-y_{\ell-m+1}\}=1$ so that
\begin{align*}
\vert  \Delta^{k}g[V_k]\vert  ^p
&\le
C(m)^p\,
\left\{A+\smed_{j=k}^{m-1}\,
\min\{1,y_{\ell+1}-y_{\ell-m+1}\}
\,\vert  \Delta^jg[y_{\ell-m+1}, \ldots ,y_{\ell-m+1+j}]\vert  ^p\right\}\\
&\le
C(m)^p\,\left\{A+A\right\}
\end{align*}
 proving inequality \rf{L-M1} in the case under consideration.

\smallskip
\textit{The second case:} $y_\ell-y_0> 1$.

We know that $1<\vert  I\vert  =y_\ell-y_0\le 2$. In this case Lemma \reff{NP-W} and \rf{X-B} tell us that
 \begin{equation}\xlabel{G-KY}
\left(\vert  \Delta^{k}g[V_k]\vert  \right)^p
=\left(\frac{1}{k!}\,\vert  G^{(k)}(\xb)\vert  \right)^p\le
2^p\,\intl_I \left(\vert  G^{(k)}(y)\vert  ^p+\vert  G^{(k+1)}(y)\vert  ^p\right)\,dy.
\end{equation}
(Note that $G^{(k+1)}\in L_p(I)$ because $G\in\LMPR$ and $k\le m-1$).

Let $J_\nu=[y_\nu,y_{\nu+m}]$, $\nu=0, \ldots ,\ell-m$.
Because $J_\nu\subset I=[y_0,y_\ell]$, we have
 \begin{equation}\xlabel{J-LN}
\vert  J_\nu\vert  \le \vert  I\vert  =y_\ell-y_0\le 2~~~\text{for every}~~~\nu=0, \ldots ,\ell-m.
\end{equation}

Let us apply Lemma \reff{K-I2} to the interval $J_\nu$, points $z_i=y_{\nu+i}$, $i=0, \ldots ,m-1$, and the function $F=G$. This lemma tells us that
 \begin{eqnarray}
\intl_{J_\nu} \,\vert  G^{(k)}(y)\vert  ^p\,dy
&\le&
C(m)^p\,
\left\{
\vphantom{\sum_{j=k}^{m-1}}
\vert  J_\nu\vert  ^{(m-k)p}\,\intl_{J_\nu}\,\vert  G^{(m)}(s)\vert  ^p\,ds\right.\nonumber\\
&&+\left.\smed_{j=k}^{m-1}
\,\vert  J_\nu\vert  ^{(j-k)p+1}\,
\,\vert  \Delta^jg[y_\nu, \ldots ,y_{\nu+j}]\vert  ^p\right\}
\nn\\
&\le&
C(m)^p\,
\left\{\intl_{J_\nu}\,\vert  G^{(m)}(s)\vert  ^p\,ds
+\smed_{j=k}^{m-1}
\,\vert  J_\nu\vert  \,
\,\vert  \Delta^jg[y_\nu, \ldots ,y_{\nu+j}]\vert  ^p\right\}.
\nn
\end{eqnarray}
 Because $\vert  J_\nu\vert  =y_{\nu+m}-y_\nu\le 2$, see \rf{J-LN}, for every $\nu=0, \ldots ,\ell-m$ the following inequality
 \begin{eqnarray}\xlabel{GK-1}
\intl_{J_\nu} \,\vert  G^{(k)}(y)\vert  ^p\,dy
&\le&
C(m)^p\,
\left\{
\vphantom{\sum_{j=k}^{m-1}}
\intl_{J_\nu}\,\vert  G^{(m)}(s)\vert  ^p\,ds\right.\nn\\
&&+\left.\smed_{j=k}^{m-1}
\,\min\{1,y_{\nu+m}-y_\nu\}
\,\vert  \Delta^jg[y_\nu, \ldots ,y_{\nu+j}]\vert  ^p\right\}
\end{eqnarray}
holds.

In the same way we prove that for every $k\le m-2$ and $\nu=0, \ldots ,\ell-m$, we have
 \begin{align}\xlabel{GK-2}
\intl_{J_\nu} \,\vert  G^{(k+1)}(y)\vert  ^p\,dy
&\le
C(m)^p\,
\left\{
\vphantom{\sum_{j=k}^{m-1}}
\intl_{J_\nu}\,\vert  G^{(m)}(s)\vert  ^p\,ds\right.\nonumber\\
&+\left.\smed_{j=k+1}^{m-1}
\,\min\{1,y_{\nu+m}-y_\nu\}
\,\vert  \Delta^jg[y_\nu, \ldots ,y_{\nu+j}]\vert  ^p\right\}.
\end{align}

Finally, we note that $I=[y_0,y_\ell]\subset \cup\{J_\nu:\nu=0, \ldots ,\ell-m\}$. This inclusion, inequalities \rf{G-KY}, \rf{GK-1} and \rf{GK-2} imply that
 \begin{eqnarray}
\vert  \Delta^{k}g[V_k]\vert  ^p
&\le&
C(m)^p\,
\left\{
\sum_{j=k}^{m-1}
\smed_{\nu=0}^{\ell-m}
\intl_{J_\nu} \vert  G^{(m)}(s)\vert  ^p\,ds
+\smed_{\nu=0}^{\ell-m}\,\smed_{j=k}^{m-1}
\,\min\{1,y_{\nu+m}-y_\nu\}
\,\vert  \Delta^jg[y_\nu, \ldots ,y_{\nu+j}]\vert  ^p\right\}\nn\\
&=&
C(m)^p\,\{A_1+A_2\}.\nn
\end{eqnarray}

Clearly,

$$
A_2=\smed_{\nu=0}^{\ell-m}\,\smed_{j=k}^{m-1}
\,\min\{1,y_{\nu+m}-y_\nu\}
\,\vert  \Delta^jg[y_\nu, \ldots ,y_{\nu+j}]\vert  ^p
=\smed_{j=k}^{m-1}\,\smed_{i=0}^{\ell-m}\,
\,\min\{1,y_{i+m}-y_i\}
\,\vert  \Delta^jg[y_i, \ldots ,y_{i+j}]\vert  ^p
$$
so that $A_2\le A$ (see definition \rf{A-D}).

Because the covering multiplicity of the family $\{J_\nu\}_{\nu=0}^{\ell-m}$ is bounded by $m+1$,
\[
A_1=\smed_{\nu=0}^{\ell-m}
\intl_{J_\nu} \vert  G^{(m)}(s)\vert  ^p\,ds
\le (m+1)\,\intl_{\R} \vert  G^{(m)}(s)\vert  ^p\,ds
=(m+1)\,\| G\| ^p_{\LMPR}\le C(m)^p\,A.
\]
 See \rf{G-A1}. Hence,
\[
\vert  \Delta^{k}g[y_0, \ldots ,y_{k}]\vert  ^p
\le C(m)^p\,\{A_1+A_2\}\le C(m)^p\,A,
\]
 and the proof of the lemma is complete.\bx

\begin{lemma}
\lbl{SLF} Let $\{x_i\}_{i=0}^n$, $n\ge m$, be a finite strictly increasing sequence in $\R$, and let $S=\{x_0, \ldots ,x_n\}$. Then for every function $g$ on $S$ and every $k\in\{0, \ldots ,m\}$ the following inequality
$$
\smed_{i=0}^{n-k}
\min\left\{1,x_{i+m}-x_{i+k-m}\right\}
\,\left|\Delta^kg[x_i,...,x_{i+k}]\right|^p\le C(m)^p\,\smed_{j=k}^{m}\,\,\smed_{i=0}^{n-j}
\min\left\{1,x_{i+m}-x_{i}\right\}
\,\left|\Delta^jg[x_i,...,x_{i+j}]\right|^p
$$
 holds. (We recall that, according to Convention \reff{AGREE-2}, $x_i=-\infty$, if $i<0$, and $x_i=+\infty$ if $i>n$).
\end{lemma}

{\it Proof.} For $k=m$ the lemma is obvious. Fix $k\in\{0, \ldots ,m-1\}$ and set
 \begin{equation}\xlabel{AL-2}
\alpha_k=\smed_{j=k}^{m}\,\,\smed_{i=0}^{n-j}
\min\left\{1,x_{i+m}-x_{i}\right\}
\,\left\lvert  \Delta^jg[x_i, \ldots ,x_{i+j}]\right\rvert  ^p
\end{equation}
 and
\[
A_k=\smed_{i=0}^{n-k}
\min\left\{1,x_{i+m}-x_{i+k-m}\right\}
\,\left\lvert  \Delta^kg[x_i, \ldots ,x_{i+k}]\right\rvert  ^p.
\]
 Clearly, $A_k\le\alpha_k+B_k$ where
\[
B_k=\,\smed_{i=0}^{n-k}
\min\left\{1,x_{i+k}-x_{i+k-m}\right\}
\,\left\lvert  \Delta^kg[x_i, \ldots ,x_{i+k}]\right\rvert  ^p.
\]

Prove that $B_k\le C(m)^p\,\alpha_k$.

We introduce a partition $\{\Ic_j:j=1,2,3\}$ of the set $\{0, \ldots ,n-k\}$ as follows: Let
\begin{align*}
\Ic_1&=\{i\in\{0, \ldots ,n-k\}: x_{i+k}-x_{i+k-m}\le 2\}, \\
\Ic_2&=\{i\in\{0, \ldots ,n-k\}: x_{i+k}-x_{i+k-m}>2~~\text{and}~~x_{i+m}-x_i>1\}
\end{align*}
 and
\[
\Ic_3=\{i\in\{0, \ldots ,n-k\}: x_{i+k}-x_{i+k-m}>2~~\text{and}~~x_{i+m}-x_i\le 1\}.
\]

Let
 \begin{equation}\xlabel{BK-J}
B_k^{(\nu)}=\,\smed_{i\in \Ic_\nu}
\min\left\{1,x_{i+k}-x_{i+k-m}\right\}
\,\left\lvert  \Delta^kg[x_i, \ldots ,x_{i+k}]\right\rvert  ^p,~~~~\nu=1,2,3
\end{equation}
   provided $\Ic_\nu\ne\emp$; otherwise, we set $B_k^{(\nu)}=0$.

Clearly, $B_k=B_k^{(1)}+B_k^{(2)}+B_k^{(3)}$. Prove that $B_k^{(\nu)}\le C(m)^p\,\alpha_k$  for every $\nu=1,2,3$.
\msk

Without loss of generality, we may assume that $\Ic_\nu\ne\emp$ for each $\nu=1,2,3$.

\smallskip
First, we estimate $B_k^{(1)}$. We note that $i+k-m\ge 0$ for each $i\in \Ic_1$; otherwise $x_{i+k-m}=-\infty$ (see
Convention \reff{AGREE-2}) which contradicts the inequality $x_{i+k}-x_{i+k-m}\le 2$.

Fix $i\in \Ic_1$. Let $z_j=x_{i+k-m+j}$, $j=0, \ldots ,m$, and let $Z=\{z_0, \ldots ,z_m\}$. Hence, $\diam Z= x_{i+k}-x_{i+k-m}\le 2$ because $i\in \Ic_1$. We apply Lemma \reff{DF-CD} to the set $Z$ and set $S=\{x_i, \ldots ,x_{i+k}\}$, and obtain the following:
 \begin{eqnarray}
\vert  \Delta^{k}g[x_i, \ldots ,x_{i+k}]\vert  ^p
&\le& C(m)^p\,
\left\{\smed_{j=k}^{m}
\,(\diam Z)^{j-k}\,\,\vert  \Delta^jg[z_0, \ldots ,z_j]\vert  \right\}^p
\nn\\
&\le&
C(m)^p\,\smed_{j=k}^{m}\,\vert  \Delta^jg[z_0, \ldots ,z_j]\vert  ^p\nn\\
&=&
C(m)^p\,\smed_{j=k}^{m}\,
\vert  \Delta^jg[x_{i+k-m}, \ldots ,x_{i+k-m+j}]\vert  ^p.\nn
\end{eqnarray}

This inequality and \rf{AL-2} enable us to estimate $B_k^{(1)}$ as follows:
 \begin{eqnarray}
B_k^{(1)}&=&\,\smed_{i\in \Ic_1}
\min\left\{1,x_{i+k}-x_{i+k-m}\right\}
\,\left\lvert  \Delta^kg[x_i, \ldots ,x_{i+k}]\right\rvert  ^p\nn\\
&\le&\,C(m)^p\,
\smed_{i\in \Ic_1}
\min\left\{1,x_{i+k}-x_{i+k-m}\right\}
\,
\smed_{j=k}^m\,\vert  \Delta^{j}g[x_{i+k-m}, \ldots ,x_{i+k-m+j}]\vert  ^p
\nn\\
&=&\,C(m)^p\,
\smed_{j=k}^m\,\smed_{i\in \Ic_1}
\min\left\{1,x_{i+k}-x_{i+k-m}\right\}
\,
\vert  \Delta^{j}g[x_{i+k-m}, \ldots ,x_{i+k-m+j}]\vert  ^p\nn\\
&\le& C(m)^p\,\alpha_k.
\nn
\end{eqnarray}

Let us estimate $B_k^{(2)}$. We know that  $x_{i+k}-x_{i+k-m}>2$ (because $i\in \Ic_2$). Hence,
\[
B_k^{(2)}=\,\smed_{i\in \Ic_2}
\,\left\lvert  \Delta^kg[x_i, \ldots ,x_{i+k}]\right\rvert  ^p.
\]
 See \rf{BK-J}. We also know that $x_{i+m}-x_i>1$ for each
$i\in \Ic_2$ so that
\[
B_k^{(2)}=\,\smed_{i\in \Ic_2}
\,\min\{1,x_{i+m}-x_i\}
\,\left\lvert  \Delta^kg[x_i, \ldots ,x_{i+k}]\right\rvert  ^p.
\]
 This equality and definition \rf{AL-2} imply that $B_k^{(2)}\le \alpha_k$.

It remains to prove that $B_k^{(3)}\le C(m)^p\,\alpha_k$. We recall that
 \begin{equation}\xlabel{I3-D}
x_{i+k}-x_{i+k-m}>2~~\text{and}~~x_{i+m}-x_i\le 1
~~\text{for  every}~~i\in\Ic_3.
\end{equation}

Let $T_i=[x_{i+k-m},x_{i+k}],$ $i\in \Ic_3$. We recall that $x_{i+k-m}=-\infty$ if $i+k-m<0$ according to our convention adopted in the formulation of the lemma. Thus, $T_i=(-\infty,x_{i+k}]$ for each $i\in \Ic_3$, $i<m-k$.
\smsk

Let $\Tc=\{T_i:~i\in \Ic_3\}$. Note that
$\vert  T\vert  >2$ for each $T\in\Tc$, see \rf{I3-D}. We also note that $i\le n-m$ provided $i\in \Ic_3$; in fact, otherwise $i+m>n$ and $x_{i+m}=+\infty$ (according to Convention \reff{AGREE-2}) which contradicts to inequality $x_{i+m}-x_i\le 1$. In particular, $i+k\le n$ for every $i\in \Ic_3$ so that $x_{i+k}<+\infty$, $i\in \Ic_3$. Thus,
each interval $T\in\Tc$ is bounded from above.

It is also clear that, given $i_0\in \Ic_3$, there are at most $2m+1$ intervals $T_i$ from $\Tc$ such that $T_{i_0}\cap T_i\ne\emp$. This property and Lemma \reff{GRAPH} imply the existence of subfamilies $\{\Tc_1, \ldots ,\Tc_\vkp\}$, $\vkp\le 2m+2$, of the family $\Tc$ such that:

\smallskip
(i) For each $j\in\{1, \ldots ,\vkp\}$ the intervals of the family $\Tc_j$ are pairwise disjoint;
\smallskip

(ii) $\Tc_{j_1}\cap\Tc_{j_2}=\emp$ for distinct $j_1,j_2\in\{0, \ldots ,\vkp\}$;
~(iii) $\Tc=\cup\{\Tc_j:j=0, \ldots ,\vkp\}$.
\smallskip

Fix $j\in\{0, \ldots ,\vkp\}$ and consider the family $\Tc_j$. This is a finite family of pairwise disjoint closed and bounded from above  intervals in $\R$. We know that $\vert  T\vert  >2$ for each  $T\in\Tc_j$. This enables us to partition $\Tc_j$ into two families, say $\Tc_j^{(1)}$ and $\Tc_j^{(2)}$, with the following properties: \textit{the distance between any two distinct intervals $T,T'\in \Tc_j^{(\nu)}$ is at least $2$.} (Here $\nu=1$ or $2$.)

For instance, we can produce the families $\Tc_j^{(1)}$ and $\Tc_j^{(2)}$ by enumerating the intervals from $\Tc_j$ in ``increasing order'', and then, by setting  $\Tc_j^{(1)}$ to be the family of intervals from $\Tc_j$ with the odd index, and $\Tc_j^{(2)}$ to be the family of intervals from $\Tc_j$ with the even index.
\smsk

These observations enable us to make the following assumption: \textit{Let $\Tc=\{T_i:~i\in \Ic_3\}$. Then}
 \begin{equation}\xlabel{DT-2}
\dist(T,T')> 2~~~~\text{for every}~~~T,T'\in\Tc,~T\ne T'.
\end{equation}

We recall that $T_i=[x_{i+k-m},x_{i+k}]$, so that   $x_i\in T_i$ for all $i\in I_3$. This and property \rf{DT-2} imply that $\vert  x_i-x_j\vert  >2$ provided $i,j\in \Ic_3$ and $i\ne j$. Hence,
 \begin{equation}\xlabel{T-IJ}
[x_i,x_i+2]\cap [x_j,x_j+2]=\emp,~~~i,j\in \Ic_3,~i\ne j.
\end{equation}

We also note that, given $i\in \Ic_3$ there exists a positive integer $\ell_i\in[m,n]$ such that
 \begin{equation}\xlabel{X-2}
x_{\ell_i}-x_i\le 2~~~~\text{but}~~~x_{\ell_i+1}-x_i>2.
\end{equation}
 This is immediate from inequality $x_{i+m}-x_i\le 1$,
$i\in\Ic_3$, see \rf{I3-D}. (In general, it may happen that $\ell_i=n$; in this case, according to Convention \reff{AGREE-2}, $x_{\ell_i+1}=+\infty$.)

Let us apply Lemma \reff{P-WQ} to the function $g$, points $y_i=x_{i+j}$, $j=0, \ldots ,n-i$, and the number $\ell=\ell_i$. This lemma tells us that the following inequality
\[
\vert  \Delta^kg[x_i, \ldots ,x_{i+k}]\vert  ^p\le  C(m)^p\,
\smed_{j=k}^{m}\,\,\smed_{\nu=i}^{\ell_i-j}
\min\left\{1,x_{\nu+m}-x_{\nu}\right\}
\,\left\lvert  \Delta^jg[x_\nu, \ldots ,x_{\nu+j}]\right\rvert  ^p
\]
 holds. This inequality and the property $x_{i+k}-x_{i+k-m}>2$, $i\in I_3$, imply that
 \begin{eqnarray}
B_k^{(3)}&=&\,\smed_{i\in \Ic_3}
\min\left\{1,x_{i+k}-x_{i+k-m}\right\}
\,\left\lvert  \Delta^kg[x_i, \ldots ,x_{i+k}]\right\rvert  ^p=
\smed_{i\in \Ic_3}
\,\left\lvert  \Delta^kg[x_i, \ldots ,x_{i+k}]\right\rvert  ^p
\nn\\
&\le& C(m)^p\,
\smed_{j=k}^{m}\,\smed_{i\in\Ic_3}
\,\smed_{\nu=i}^{\ell_i-j}
\min\left\{1,x_{\nu+m}-x_{\nu}\right\}
\,\left\lvert  \Delta^jg[x_\nu, \ldots ,x_{\nu+j}]\right\rvert  ^p=
C(m)^p\,
\smed_{j=k}^{m}\,Y_j\,.\nn
\end{eqnarray}

Note that, given $j\in\{0, \ldots ,m\}$ and $i\in\Ic_3$, the points $x_i, \ldots , x_{\ell_i}\in[x_i,x_i+2]$, see \rf{X-2}. Thanks to this property and \rf{T-IJ}, for every $i',i''\in\Ic_3$, $i'\ne i''$, the families of indexes
$\{i', \ldots ,\ell_{i'}\}$ and $\{i'', \ldots ,\ell_{i''}\}$ are disjoint. Hence,
\begin{align*}
Y_j&=\,\smed_{i\in\Ic_3}
\,\smed_{\nu=i}^{\ell_i-j}
\min\left\{1,x_{\nu+m}-x_{\nu}\right\}
\,\left\lvert  \Delta^jg[x_\nu, \ldots ,x_{\nu+j}]\right\rvert  ^p\\
&\le
\,\smed_{i=0}^{n-j}
\min\left\{1,x_{i+m}-x_{i}\right\}
\,\left\lvert  \Delta^jg[x_i, \ldots ,x_{i+j}]\right\rvert  ^p.
\end{align*}
 From this inequality and \rf{AL-2}, we have
\[
B_k^{(3)}\le C(m)^p\,
\smed_{j=k}^{m}\,\,\smed_{i=0}^{n-j}
\min\left\{1,x_{i+m}-x_{i}\right\}
\,\left\lvert  \Delta^jg[x_i, \ldots ,x_{i+j}]\right\rvert  ^p=
C(m)^p\,\alpha_k.
\]

Thus, $B_k^{(\nu)}\le C(m)^p\,\alpha_k$, $\nu=1,2,3$,
so that
$B_k=B_k^{(1)}+B_k^{(2)}+B_k^{(3)}\le C(m)^p\,\alpha_k$.

Hence, $A_k\le \alpha_k+B_k\le  C(m)^p\,\alpha_k$, and the proof of the lemma is complete.\bx


\SECT{4. The sufficiency part of the variational extension criterion: proofs.}{4}
\addtocontents{toc}{4. The sufficiency part of the variational extension criterion: proofs. \hfill\thepage\par \VST}

\indent\par In this section we proof the sufficiency part of Theorem \reff{W-VAR-IN}. Let $p\in(1,\infty]$, $m\in\N$, and let $E$ be a closed subset of $\R$ with $\#E\ge m+1$. Let $f$ be a function on $E$ such that
 \begin{equation}\xlabel{LMBD}
\lambda=\NWMP(f:E)<\infty.
\end{equation}
 See \rf{N-WRS-IN} and \rf{N-WINF}. These definitions enable us to make the following assumption.
\begin{assumption}
\lbl{ASMP-FE} Let $1<p<\infty$. Then for every integer $n\ge m$ and every strictly increasing sequence $\{x_i\}_{i=0}^n$ in $E$ the following inequality
\[
\smed_{k=0}^{m}\,\,\smed_{i=0}^{n-k}
\min\left\{1,x_{i+m}-x_{i}\right\}
\,\left\lvert  \Delta^kf[x_i, \ldots ,x_{i+k}]\right\rvert  ^p
\le \lambda^p
\]
 holds.

Let $p=\infty$. Then for every $k\in\{0, \ldots ,m\}$ and every $(k+1)$-point subset $\{y_0, \ldots ,y_k\}\subset E$ we have
 \begin{equation}\xlabel{A-INF}
\vert  \Delta^kf[y_0, \ldots ,y_k]\vert  \le \lambda.
\end{equation}
\end{assumption}
\smsk

Our aim is to prove that $ f\in\WMPR\vert  _E$ and $\| f\| _{\WMPR\vert  _E}\le C(m)\,\lambda$.
\smallskip

Let us construct an almost optimal $\WMPR$-extension of the function $f$.

Because $E$ is a closed subset of $\R$, the complement of $E$, the set $\R\setminus E$, can be represented as a union of a certain finite or countable family
$\Jc_E=\{J_k=(a_k,b_k): k\in \Kc\}$ of pairwise disjoint open intervals (bounded or unbounded). Thus, $a_k,b_k\in E\cup\{\pm\infty\}$ for all $k\in\Kc$,
\[
\R\setminus E= \mcup\{J_k=(a_k,b_k): k\in\Kc\}~~~\text{and}~~~J_{k'}\mcap\, J_{k''}=\emp~~
\text{for every}~~k',k''\in \Kc, k'\ne k''.
\]

Next, we introduce a (perhaps empty) subfamily $\Gc_E$ of $\Jc_E$ defined by
\[
\Gc_E=\{J\in\Jc_E:\vert  J\vert  >4\}.
\]

Given a \textit{bounded} interval $J=(a_J,b_J)\in\Gc_E$, we put
 \begin{equation}\xlabel{NJ-12}
n_J=\lfloor\vert  J\vert  /2\rfloor
\end{equation}
 where $\lfloor \cdot\rfloor$ denotes the greatest integer function. Let $\ell(J)=\vert  J\vert  /n_J$; then $2\le\ell(J)\le 3$.

Let us associate to the interval $J$ points $Y^{(J)}_n\in J$ defined by
 \begin{equation}\xlabel{YNJ}
Y^{(J)}_n=a_J+\ell(J)\cdot n,~~~~n=1, \ldots ,n_J-1.
\end{equation}
 We also set $Y^{(J)}_0=a_J$ and $Y^{(J)}_{n_J}=b_J$ and put
\[
S_J=\{Y^{(J)}_1, \ldots ,Y^{(J)}_{n_J-1}\}~~~
\text{for every bounded interval}~~~J\in \Gc_E.
\]

Thus, the points $Y^{(J)}_0, \ldots ,Y^{(J)}_{n_J}$ divide $J=(a_J,b_J)$ in $n_J$ subintervals $(Y^{(J)}_n,Y^{(J)}_{n+1})$, $n=0, \ldots ,n_J-1$, of equal length ($=\ell(J)$). We know that
 \begin{equation}\xlabel{LJ-23}
2\le \ell(J)=Y^{(J)}_{n+1}-Y^{(J)}_n\le 3~~~~\text{for every}~~~~n=0, \ldots ,n_J-1.
\end{equation}

Let $J=(a_J,b_J)$ be an \textit{unbounded} interval. (Clearly, $J\in\Gc_E$). In this case we set
\[
Y^{(J)}_n=a_J+2\cdot n,~~~~n\in \N,
\]
 provided $J=(a_J,+\infty)$, and
 \begin{equation}\xlabel{LJ-MI}
Y^{(J)}_n=b_J-2\cdot n,~~n\in \N,
~~~~~\text{provided}~~~J=(-\infty,b_J).
\end{equation}

In other words, we divide each unbounded interval $J$ in subintervals of length $2$. Finally, we set
\[
S_J=\{Y^{(J)}_n:n\in\N\}~~~\text{for every unbounded interval}~~~J\in \Gc_E.
\]

Let
 \begin{equation}\xlabel{G-WE}
G=\bigcup_{J\in\Gc_E} S_J
\end{equation}
 whenever $\Gc_E\ne\emp$, and $G=\emp$  otherwise. Clearly,
 \begin{equation}\xlabel{D-EG}
\dist(E,G)\ge 2.
\end{equation}
 (Recall that $\dist(A,\emp)=+\infty$ provided $A\ne\emp$, so that \rf{D-EG} includes the case of $G=\emp$ as well).
We also note that \rf{LJ-23}--\rf{LJ-MI} imply that
 \begin{equation}\xlabel{SP-YJ}
\vert  Y-Y'\vert  \ge 2~~~~\text{for every}~~~~Y,Y'\in G, Y\ne Y'.
\end{equation}

Let
 \begin{equation}\xlabel{EW-G}
\tE=E\cup G.
\end{equation}

Note that from \rf{YNJ}--\rf{G-WE} and \rf{EW-G}, we have
 \begin{equation}\xlabel{EW-NBH}
\dist(x,\tE)\le 2~~~\text{for every}~~~x\in\R.
\end{equation}

By $\tf:\tE\to \R$ we denote \textit{the extension of $f$ from $E$ to $\tE$ by zero}; thus,
 \begin{equation}\xlabel{TF-D2}
\tf\vert  _E=f~~~~~~\text{and}~~~~~~\tf\vert  _G\equiv 0.
 \end{equation}

Let us show that for every $p\in(1,\infty]$
\[
\tf\in \WMPR\vert  _{\tE}~~~\text{and}~~~\| \tf\| _{\WMPR\vert  _{\tE}}\le C(m)\,\lambda
\]
 provided Assumption \reff{ASMP-FE} holds.

Our proof of this statement relies on a series of auxiliary lemmas.
\begin{lemma}
 \lbl{SM-E} Let $p\in(1,\infty]$ and let $k\in\{0, \ldots ,m-1\}$. Then for every strictly increasing sequence $\{y_i\}_{i=0}^k$ in $E$ the following inequality
 \begin{equation}\xlabel{F-LS1}
\vert  \Delta^kf[y_0, \ldots ,y_k]\vert  \le C(m)\,\lambda
\end{equation}
 holds.
\end{lemma}

{\it Proof.} The lemma is obvious for $p=\infty$. See \rf{A-INF}. Let $1<p<\infty$. We know that $E$ contains at least $m+1$ distinct points. Therefore, there exists a strictly increasing sequence $\{x_\nu\}_{\nu=0}^m$ in $E$ containing the points $y_0, \ldots ,y_k$. Thus, the set $Y=\{y_0, \ldots ,y_k\}$ is  a $(k+1)$-point subset of the set $X=\{x_0, \ldots ,x_m\}\subset E$ such that $x_0<\cdots<x_m$.

Property $(\bigstar 3)$ of Sections 2 (with $n=m$) tells us that there exist non-negative numbers $\alpha_\nu$, $\nu=0, \ldots ,m-k$, such that $\alpha_0+\cdots+\alpha_{m-k}=1$, and
\[
\Delta^kf[Y]=\smed_{\nu=0}^{m-k}\,\alpha_\nu\,
\Delta^kf[x_\nu, \ldots ,x_{\nu+k}].
\]
 Hence,
\[
\vert  \Delta^kf[Y]\vert  ^p\le (m+1)^p
\smed_{\nu=0}^{m-k}\,\vert  \Delta^kf[x_\nu, \ldots ,x_{\nu+k}]\vert  ^p.
\]

We note that for every $\nu\in\{0, \ldots ,m-k\}$ either $\nu+m>m$ or $v+k-m<0$ (because $0\le k<m$). Therefore, according to Convention \reff{AGREE-2}, either $x_{\nu+m}=+\infty$ or $x_{\nu+k-m}=-\infty$ proving that
\[
\min\{1,x_{\nu+m}-x_{\nu+k-m}\}=1~~~\text{for all}~~~\nu=0, \ldots ,m-k.
\]

Thus,
\[
\vert  \Delta^kf[Y]\vert  ^p\le (m+1)^p
\smed_{\nu=0}^{m-k}\,\min\{1,x_{\nu+m}-x_{\nu+k-m}\}\,
\vert  \Delta^kf[x_\nu, \ldots ,x_{\nu+k}]\vert  ^p.
\]

This inequality and Lemma \reff{SLF} (with $n=m$) imply that
\[
\vert  \Delta^kf[Y]\vert  ^p\le C(m)^p\,
\smed_{j=k}^{m}\smed_{i=0}^{m-j}
\,\min\{1,x_{i+m}-x_{i}\}\,\vert  \Delta^jf[x_i, \ldots ,x_{i+j}]\vert  ^p.
\]

Finally, applying Assumption \reff{ASMP-FE} to the right hand side of this inequality, we get the required inequality \rf{F-LS1} proving the lemma.\bx

\begin{lemma}
\lbl{LM-EW} For every finite strictly increa\-sing sequence of points $\{y_i\}_{i=0}^n\subset \tE$, $n\ge m$, the following inequality
\[
\smed_{i=0}^{n-m}\,
(y_{i+m}-y_{i})\left\lvert  \Delta^m\tf[y_i, \ldots ,y_{i+m}]\right\rvert  ^p
\le\,C(m)^p\,\lambda^p
\]
 holds provided $1<p<\infty$. If $p=\infty$, then for every $(m+1)$ point set $\{y_0, \ldots ,y_m\}\subset \tE$ we have
\[
\vert  \Delta^m\tf[y_0, \ldots ,y_{m}]\vert  \le\,C(m)\,\lambda.
\]
\end{lemma}

{\it Proof.} We may assume that $m\ge 1$ (for $m=0$ the lemma is trivial). First we prove the lemma for $p\in(1,\infty)$. Let
\[
A=\smed_{i=0}^{n-m}\,
(y_{i+m}-y_{i})\left\lvert  \Delta^m\tf[y_i, \ldots ,y_{i+m}]\right\rvert  ^p
\]
 and let
\begin{equation}\xlabel{IGE-12}
I_1=\{i:0\le i\le n-m, y_{i+m}-y_{i}<2\},~
I_2=\{i:0\le i\le n-m, y_{i+m}-y_{i}\ge 2\}.
\end{equation}
 Given $j\in\{1,2\}$, let
\[
A_j=\smed_{i\in I_j}\,
(y_{i+m}-y_{i})\left\lvert  \Delta^m\tf[y_i, \ldots ,y_{i+m}]\right\rvert  ^p
\]
provided $I_j\ne\emp$, and let $A_j=0$ otherwise. Clearly, $A=A_1+A_2$.

Prove that
 \begin{equation}\xlabel{L1}
A_1\le 2\,\lambda^p.
\end{equation}

This inequality is trivial if $I_1=\emp$. Let us assume that  $I_1\ne\emp$. Then, thanks to \rf{SP-YJ} and definition of $I_1$ (see \rf{IGE-12}), $\{y_i:i\in I_1\}\subset E$. Furthermore, we know that $y_{i+m}-y_{i}<2$ for every $i\in I_1$, and $\dist(E,G)\ge 2$ (see \rf{D-EG}). These inequalities imply that $y_i, \ldots ,y_{i+m}\in E$ for every $i\in I_{1}$.

We also recall that $\tf\vert  _E=f$. Therefore,
\begin{eqnarray}
A_1&=&
\smed_{i\in I_1}\,
(y_{i+m}-y_{i})\left\lvert  \Delta^m\tf[y_i, \ldots ,y_{i+m}]\right\rvert  ^p
=
\smed_{i\in I_{1}}\,
(y_{i+m}-y_{i})\left\lvert  \Delta^mf[y_i, \ldots ,y_{i+m}]\right\rvert  ^p
\nn\\
&\le& 2\,
\smed_{i\in I_{1}}\,\min\{1,y_{i+m}-y_{i}\}\,
\left\lvert  \Delta^mf[y_i, \ldots ,y_{i+m}]\right\rvert  ^p.\nn
\end{eqnarray}
 Hence,
\[
A_1\le
2\,\smed_{i=0}^{n-m}\,
\min\{1,x_{i+m}-x_{i}\}\,
\left\lvert  \Delta^mf[x_i, \ldots ,x_{i+m}]\right\rvert  ^p,
\]
 which together with Assumption \reff{ASMP-FE} implies the required inequality \rf{L1}.
\msk

Prove that
 \begin{equation}\xlabel{L-A2}
A_2\le C(m)^p\,\lambda^p.
\end{equation}

We may assume that the set $I_2$ determined in \rf{IGE-12} is not empty; otherwise, $A_2=0$.

Let $i\in I_2$ and let $T_i=[y_{i-m},y_{i+2m}]$. (We recall our Convention \reff{AGREE-2} concerning the values of $y_j$ for $i\notin \{0, \ldots ,n\}$.) Let
$\Tc=\{T_i: i\in I_2\}$.
\smsk

Note that for each interval $T_{i_0}\in\Tc$, there exist at most $6m+1$ intervals $T_i\in\Tc$ such that $T_i\cap T_{i_0}\ne\emp$. Lemma \reff{GRAPH} tells us that we can partition $\Tc$ in at most $\ell\le 6m+2$ subfamilies $\{\Tc_1, \ldots ,\Tc_{\ell}\}$ each consisting of pairwise disjoint intervals. This enables us, without loss of generality, to assume that \textit{the family $\Tc$ itself consists of pairwise disjoint intervals}, i.e.,~\begin{equation}\xlabel{T-ASM}
T_i\cap T_j\ne\emp~~~\text{for every}~~~i,j\in I_2,
~i\ne j.
\end{equation}
 In particular, this property implies that
$\vert  i-j\vert  >3m$ for every $i,j\in I_2,~i\ne j$.
\smsk

Fix $i\in I_2$ and apply Lemma \reff{SM-DF} to points $\{y_i, \ldots ,y_{i+m}\}$ and the function $\tf$. (Recall that $y_{i+m}-y_i\ge 2$ for all $i\in I_2$.) This lemma tells us that there exist $k_i\in \{0, \ldots ,m-1\}$ and $\alpha_i\in\{i, \ldots ,i+m-k_i\}$ such that $y_{\alpha_i+k_i}-y_{\alpha_i}\le 1$ and
 \begin{equation}\xlabel{DM-1}
\vert  \Delta^m\tf[y_i, \ldots ,y_{i+m}]\vert  \le\,2^m \,\vert  \Delta^{k_i}\tf[y_{\alpha_i}, \ldots ,y_{\alpha_i+k_i}]\vert  /
(y_{i+m}-y_i).
\end{equation}
 Furthermore, either
\begin{equation}\xlabel{AKI}
\alpha_i+k_i+1\le i+m~\text{and}~
y_{\alpha_i+k_i+1}-y_{\alpha_i}\ge 1,~\text{or}~
\alpha_i\ge i+1~\text{and}~
y_{\alpha_i+k_i}-y_{\alpha_i-1}\ge 1.
\end{equation}

Because $y_{i+m}-y_i\ge 2$ and $p>1$, from inequality \rf{DM-1} it follows that
\begin{align*}
(y_{i+m}-y_i)\vert  \Delta^m\tf[y_i, \ldots ,y_{i+m}]\vert  ^p
&\le\,2^{mp} \,\vert  \Delta^{k_i}\tf[y_{\alpha_i}, \ldots ,y_{\alpha_i+k_i}]\vert  ^p/
(y_{i+m}-y_i)^{p-1}\\
&\le \,2^{mp} \,\vert  \Delta^{k_i}\tf[y_{\alpha_i}, \ldots ,y_{\alpha_i+k_i}]\vert  ^p.
\end{align*}
 Moreover, because $y_{\alpha_i+k_i}-y_{\alpha_i}\le 1$ and $\dist(E,G)\ge 2$ (see \rf{D-EG}), either the set $\{y_{\alpha_i}, \ldots ,y_{\alpha_i+k_i}\}\subset E$ or
$\{y_{\alpha_i}, \ldots ,y_{\alpha_i+k_i}\}\subset G$.

We know that $\tf\equiv 0$ on $G$. Therefore,
 \begin{equation}\xlabel{A-2MP}
A_2=\smed_{i\in I_2}
(y_{i+m}-y_{i})\left\lvert  \Delta^m\tf[y_i, \ldots ,y_{i+m}]\right\rvert  ^p
\le 2^{mp}\smed_{i\in I_{2,E}}
\vert  \Delta^{k_i}f[y_{\alpha_i}, \ldots ,y_{\alpha_i+k_i}]\vert  ^p
\end{equation}
 where
 \begin{equation}\xlabel{I-2E}
I_{2,E}=\{i\in I_2:y_{\alpha_i}, \ldots ,y_{\alpha_i+k_i}\in E\}.
\end{equation}

We recall that $k_i\in\{0, \ldots ,m-1\}$, $\alpha_i\in\{0, \ldots ,m-k_i\}$ for each $i\in I_{2,E}$ so that
\[
y_{\alpha_i+m}-y_{\alpha_i+k_i-m}\ge
\max\{y_{\alpha_i+k_i+1}-y_{\alpha_i},
y_{\alpha_i+k_i}-y_{\alpha_i-1}\}.
\]
 This inequality and property \rf{AKI} imply that
 \begin{equation}\xlabel{Y-G1}
y_{\alpha_i+m}-y_{\alpha_i+k_i-m}\ge 1 ~~~\text{for every}~~~i\in I_{2,E}.
\end{equation}

Let
 \begin{equation}\xlabel{H-U2-S5}
H=\bigcup_{i\in\, I_{2,E}}
\{y_{\alpha_i}, \ldots ,y_{\alpha_i+k_i}\}~~~~\text{and}
~~~~\vkp=\#H.
\end{equation}

We know that $H\subset E$, see \rf{I-2E}. In turn, assumption \rf{T-ASM} implies the following property:
\[
\{y_{\alpha_i}, \ldots ,y_{\alpha_i+k_i}\}\capsm
\{y_{\alpha_j}, \ldots ,y_{\alpha_j+k_j}\}=\emp~~~\text{for all}~~~i,j\in I_{2,E},~i\ne j.
\]
 In particular, $\{\alpha_i\}_{i\in I_{2,E}}$ is a strictly increasing sequence, and $\#I_{2,E}\le\#H=\vkp$.
\msk

Consider two cases.

\smallskip
\textit{The first case: $\vkp\ge m+1$}. Since $\{y_{i}\}_{i=0}^n$ is a strictly increasing sequence, and $H\subset E$, we can consider the set $H$
as a strictly increasing subsequence of the sequence $\{y_{i}\}_{i=0}^n$ whose elements lie in $E$. This enables us to represent $H$ in the form $H=\{x_0, \ldots ,x_\vkp\}$ where $x_0<\cdots<x_\vkp$.

Furthermore, we know that for each  $i\in I_{2,E}$ there exists a unique $\nu_i\in\{0, \ldots ,\vkp\}$ such that $y_{\alpha_i}=x_{\nu_i}$. The sequence $\{\nu_i\}_{i\in I_{2,E}}$ is a strictly increasing sequence such that
$y_{\alpha_i+j}=x_{\nu_i+j}$ for every $j=0, \ldots ,k_i$. See \rf{I-2E}. We also note that, according to Convention \reff{AGREE-2}, $y_{\alpha_i-j}=-\infty$ if $\alpha-j<0$ and $x_{v_i+j}=+\infty$ if $v_i+j>\kappa$. Because $\{x_\nu\}_{\nu=0}^\vkp$ is a strictly increasing subsequence of $\{y_{i}\}_{i=0}^n$, we have
 \begin{equation}\xlabel{XIV}
x_{\nu_i+j}\ge y_{\alpha_i+j}~~~\text{and}~~~
x_{\nu_i-j}\le y_{\alpha_i-j}~~~\text{for every}~~~j=0, \ldots ,m.
\end{equation}

Hence,
\[
\Delta^{k_i}f[y_{\alpha_i}, \ldots ,y_{\alpha_i+k_i}]=
\Delta^{k_i}f[x_{\nu_i}, \ldots ,x_{\nu_i+k_i}]~~~\text{for all}~~~i\in I_{2,E}.
\]

From \rf{Y-G1} and \rf{XIV} we have
$x_{\nu_i+m}-x_{\nu_i+k_i-m}\ge y_{\alpha_i+m}-y_{\alpha_i+k_i-m}\ge 1$,
so that
\[
\vert  \Delta^{k_i}f[y_{\alpha_i}, \ldots ,y_{\alpha_i+k_i}]\vert  ^p=
\min\{1,x_{\nu_i+m}-x_{\nu_i+k_i-m}\}
\vert  \Delta^{k_i}f[x_{\nu_i}, \ldots ,x_{\nu_i+k_i}]\vert  ^p, ~~~i\in I_{2,E}.
\]
 This equality and \rf{A-2MP} imply the following estimate of $A_2$:
\[
A_2\le \,2^{mp} \,\smed_{i\in I_{2,E}}\,
\min\{1,x_{\nu_i+m}-x_{\nu_i+k_i-m}\}\,
\vert  \Delta^{k_i}f[x_{\nu_i}, \ldots ,x_{\nu_i+k_i}]\vert  ^p.
\]
    Hence,
 \begin{equation}\xlabel{A2-LT}
A_2\le \,2^{mp} \,\smed_{k=0}^m\,\smed_{\nu=0}^{\tn-k}
\min\{1,x_{\nu+m}-x_{\nu+k-m}\}
\vert  \Delta^{k}f[x_{\nu}, \ldots ,x_{\nu+k}]\vert  ^p
\end{equation}
 where $\tn=\vkp-1$. We know that in the case under consideration $\tn\ge m$.

Lemma \reff{SLF} tells us that for every $k\in\{0, \ldots ,m\}$ the following inequality
$$
\smed_{\nu=0}^{\tn-k}
\min\left\{1,x_{\nu+m}-x_{\nu+k-m}\right\}
\,\left|\Delta^kf[x_\nu,...,x_{\nu+k}]\right|^p\le C(m)^p\,\smed_{j=k}^{m}\,\smed_{i=0}^{\tn-j}
\min\left\{1,x_{i+m}-x_{i}\right\}
\,\left|\Delta^jf[x_i,...,x_{i+j}]\right|^p
$$
 holds. On the other hand, thanks to Assumption \reff{ASMP-FE},
 \begin{equation}\xlabel{L2}
\smed_{j=k}^{m}\,\,\smed_{i=0}^{\tn-k}
\min\left\{1,x_{i+m}-x_{i}\right\}
\,\left\lvert  \Delta^jf[x_i, \ldots ,x_{i+k}]\right\rvert  ^p
\le \lambda^p
\end{equation}
 proving that
\[
\smed_{\nu=0}^{\tn-k}
\min\left\{1,x_{\nu+m}-x_{\nu+k-m}\right\}
\,\left\lvert  \Delta^kf[x_\nu, \ldots ,x_{\nu+k}]\right\rvert  ^p
\le C(m)^p\,\lambda^p
\]
for every $k\in\{0, \ldots ,m\}$. This inequality and \rf{A2-LT} imply the required inequality \rf{L-A2}.
\smallskip

\textit{The second case: $\vkp\le m$}. Fix $i\in I_{2,E}$. From \rf{H-U2-S5} we have \[0<k_i\le \#H-1=\vkp-1\le m-1.\] Lemma \reff{SM-E} tells us that
$
\vert  \Delta^{k_i}f[y_{\alpha_i}, \ldots ,y_{\alpha_i+k_i}]\vert  \le  C(m)\,\lambda
$
 which together with \rf{A-2MP} implies that
 $
A_2\le 2^{mp}\,C(m)^p\,\lambda^p\cdot \#I_{2,E}.
$
 But $\#I_{2,E}\le \#H=\vkp\le m$ proving \rf{L-A2} in the second case.
\smallskip

Thus, inequality \rf{L-A2} holds. This inequality together with \rf{L1} implies that
\[
A=A_1+A_2\le 2\lambda^p + C(m)^p\,\lambda^p=(2+C(m)^p)\,\lambda^p,
\]
 proving the lemma for $p\in(1,\infty)$.
\smsk

Prove the lemma for $p=\infty$. Lemma \reff{SM-DF} produces a certain $k\in\{0,\ldots,m-1\}$ and $i\in\{0, \ldots ,m-k\}$ such that $y_{i+k}-y_i\le 1$ and
\[
\vert  \Delta^m\tf[y_0, \ldots ,y_m]\vert  \le\,
2^m \,\vert  \Delta^k\tf[y_i, \ldots ,y_{i+k}]\vert  .
\]

We know that $\dist(E,G)\ge 2$ (see \rf{D-EG}). But
$\tY=\{y_i, \ldots ,y_{i+k}\}\subset \tE=E\cup G$ and $y_{i+k}-y_i\le 1$ proving that either $\tY\subset G$ or $\tY\subset E$. Clearly, if $\tY\subset G$ then $\Delta^k\tf[\tY]=0$ (because $\tf\vert  _G\equiv 0$), and therefore $\Delta^m\tf[y_0, \ldots ,y_m]=0$.

In turn, if $\tY\subset E$, then $\tf\vert  _{\tY}=f\vert  _{\tY}$ so that, thanks to \rf{A-INF},
 \begin{equation}\xlabel{PI-R}
\vert  \Delta^m\tf[y_0, \ldots ,y_m]\vert  \le\,
2^m \,\vert  \Delta^kf[y_i, \ldots ,y_{i+k}]\vert  \le 2^m \,\lambda\,.
\end{equation}

The proof of the lemma is complete.\bx
\msk

Lemma \reff{LM-EW} shows that the quantity  $\NMP(f:\tE)$ (see \rf{L-INF1} and \rf{NMP}) is bounded by $C(m)\,\lambda$ provided Assumption \reff{ASMP-FE} holds for $f$. In other words, Lemma \reff{LM-EW} proves that
\[
\NMP(f:\tE)\le C(m)\,\NWMP(f:E)~~~\text{for every}~~~p\in(1,\infty].~~~~~~\text{(See \rf{LMBD})}.
\]

Now, Theorem \reff{MAIN-TH} $(1<p<\infty)$ and equivalence \rf{T-R1} $(p=\infty)$ tell us that there exists a function $F\in\LMPR$ with $\| F\| _{\LMPR}\le C(m)\,\NMP(f:\tE)$ such that $F\vert_{\tE}=\tf$. Thus,
 \begin{equation}\xlabel{N-WMP}
\| F\| _{\LMPR}\le C(m)\,\lambda=C(m)\,\NWMP(f:E).
\end{equation}
 Since $\tf\vert  _E=f$, the function $F$ is an extension of $f$ from $E$ to all of $\R$, i.e.,~$F\vert_E=f$. We denote the extension $F$ by
 \begin{equation}\xlabel{EXT-W}
F=\EXT_E(f:\WMPR).
\end{equation}

Thus,
 \begin{equation}\xlabel{EXT-WL}
\EXT_E(f:\WMPR)=\WOPTE(\tf)
\end{equation}
 where $\WOP$ is the Whitney extension operator \rf{EXT-L}. \smsk

\smallskip
Let $\ve=3$ and let $A$ be a \textit{maximal $\ve$-net in $\tE$}; thus $A$ is a subset of $\tE$ having the following properties:
\smsk
\medskip

$({\bigstar}1)$~$\vert  a-a'\vert  \ge 3$ for every $a,a'\in A$, $a\ne a'$;
\smallskip

$({\bigstar}2)$~$\dist(y,A)<3$ for every $y\in\tE$.

\medskip
Now, \rf{EW-NBH} and $({\bigstar}2)$ imply that
 \begin{equation}\xlabel{ES-R}
\dist(x,A)\le 5~~~\text{for every}~~~x\in\R.
\end{equation}

This inequality and property $({\bigstar}1)$ enable us to represent the set $A$ as a certain bi-infinite strongly increasing sequence $\{a_i\}_{i=-\infty}^{+\infty}$ in $\tE$. Thus,
 \begin{equation}\xlabel{A-SQN}
A=\{a_i\}_{i=-\infty}^{+\infty}~~~\text{where}~~~a_i\in\tE~~~\text{and}~~~a_i<a_{i+1}~~~\text{for all}~~~i\in\Z.
\end{equation}
 Furthermore, thanks to  $({\bigstar}1)$ and \rf{ES-R},
 \begin{equation}\xlabel{AI-INF}
a_{i}\to-\infty~~~\text{as}~~~i\to-\infty,~~~\text{and}~~~
a_{i}\to+\infty~~~\text{as}~~~i\to+\infty,
\end{equation}
   and, thanks to the property $({\bigstar}1)$ and \rf{ES-R},
 \begin{equation}\xlabel{A-MD}
3\le a_{i+1}-a_i\le 10~~~\text{for every}~~~i\in\Z.
\end{equation}
\begin{lemma}
\lbl{LP-F-G} For every $p\in(1,\infty)$ and every function $g\in\LMPR$ such that $g\vert_G\equiv 0$ (see \rf{G-WE}) the following inequality
 \begin{equation}\xlabel{G-LP}
\| g\| _{\LPR}^p\le C(m)^p\,\left\{\| g\| _{\LMPR}^p+
\smed_{i=-\infty}^{i=+\infty}\,\vert  g(a_i)\vert  ^p\right\}
\end{equation}
 holds. Furthermore, if $g\in\LMIR$ and $g\vert_G\equiv 0$, then
 \begin{equation}\xlabel{G-LI1}
\| g\| _{\LIR}\le C(m)\,\left(\| g\| _{\LMIR}+
\sup\,\{\vert  g(a_i)\vert  :~i\in \Z\}\right).
\end{equation}
\end{lemma}

{\it Proof.} Fix $N>m$ and prove that
 \begin{equation}\xlabel{INT-1}
I(N)=\intl_{a_{-N}}^{a_N}\,\vert  g(t)\vert  ^p\,dt\le C(m)^p\,\left\{\| g\| _{\LMPR}^p+
\smed_{i=-\infty}^{i=+\infty}\,\vert  g(a_i)\vert  ^p\right\}.
\end{equation}

In the proof of this inequality we use some ideas of the
work~\cite[p. 451]{Es}.

Let $i\in\Z$, $-N\le i\le N-m$, and let $T_i=[a_i,a_{i+m}]$. We apply Lemma \reff{K-I2} to the interval $T_i$ and the function $g$ taking $k=0$ and $z_j=a_{i+j}$, $j=0, \ldots ,m-1$. This lemma tells us that
\[
\intl_{T_i}\vert  g(t)\vert  ^p\,dt\le C(m)^p\,
\left\{\vert  T_i\vert  ^{mp}\,\intl_{T_i}\,\vert  g^{(m)}(s)\vert  ^p\,ds
+\smed_{j=0}^{m-1}
\,\vert  T_i\vert  ^{jp+1}\,\,\vert  \Delta^jg[a_i, \ldots ,a_{i+j}]\vert  ^p\right\}.
\]

We note that $\vert  a_\nu-a_\mu\vert  \ge 3$ for every $\nu,\mu\in\{i, \ldots ,i+m-1\}$, $\nu\ne\mu$, see \rf{A-MD}, so that, thanks to \rf{D-PT1},
 \begin{equation}\xlabel{G-PIL}
\vert  \Delta^jg[a_i, \ldots ,a_{i+j}]\vert  \le \smed_{j=0}^{m-1}\,\,\vert  g(a_{i+j})\vert  .
\end{equation}
 Because $a_{i+1}-a_i\le 10$ for every $i\in\Z$ (see again \rf{A-MD}), we have
\[
\vert  T_i\vert  =a_{i+m}-a_i\le 10m.
\]
 Hence,
\[
\intl_{T_i}\vert  g(t)\vert  ^p\,dt\le C(m)^p\,
\left\{\intl_{T_i}\,\vert  g^{(m)}(s)\vert  ^p\,ds
+\smed_{j=0}^{m-1}\,\vert  g(a_{i+j})\vert  ^p\right\}
\]
for every $-N\le i\le N-m$. This inequality implies the following estimate of the quantity $I(N)$ (see \rf{INT-1}):
\[
I(N)\le\smed_{i=-N}^{N-m}
\intl_{T_i}\vert  g(t)\vert  ^pdt
\le
C(m)^p
\left\{\smed_{i=-N}^{N-m}
\intl_{T_i}\vert  g^{(m)}(s)\vert  ^pds
+\smed_{i=-N}^{N-m}
\smed_{j=0}^{m-1}\vert  g(a_{i+j})\vert  ^p\right\}.
\]

Note that, given $i_0\in\Z$, $-N\le i_0\le N$, there exist at most $2m+1$ intervals from the family of intervals $\{T_i: i=-N, \ldots ,N-m\}$ which have common points with the interval $T_{i_0}$. Therefore,
$$
I(N)\le C(m)^p\,\left\{\,
\intl_{a_{-N}}^{a_N}\,\vert  g^{(m)}(s)\vert  ^p\,ds+
\smed_{i=-N}^{N}\,\vert  g(a_{i})\vert  ^p\right\}\le C(m)^p\,\left\{\,
\intl_{\R}\,\vert  g^{(m)}(s)\vert  ^p\,ds+
\smed_{i=-\infty}^{+\infty}\,\vert  g(a_{i})\vert  ^p\right\}
$$
which yields \rf{INT-1}. But, thanks to \rf{AI-INF}, $I(N)\to \| g\| _{\LPR}^p$ as $N\to\infty$ which implies inequality \rf{G-LP} and proves the lemma for $p\in(1,\infty)$.

For the case $p=\infty$ the required inequality \rf{G-LI1} is immediate from Lemma \reff{K-I1} (with $p=1$) and inequality \rf{G-PIL}. We leave the details to the interested reader.

The proof of the lemma is complete.\bx

\begin{lemma}
\lbl{LP-F} $\| F\| _{\LPR}\le C(m)\,\lambda$ for every $p\in(1,\infty]$.
\end{lemma}

{\it Proof.} Let $p=\infty$. Then, thanks to assumption \rf{A-INF} (with $k=0$), $\vert  f(x)\vert  \le\lambda$ for every $x\in E$. We know that $F\vert_G=\tf\vert  _G\equiv 0$, and $F\vert_E=\tf\vert  _E=f$, so that
\[
\sup\{\vert  F(a_i)\vert  :~i\in\Z\}=
\sup\{\vert  f(a_i)\vert  :~i\in\Z,\,a_i\in E\}\le\lambda.
\]

Now, Lemma \reff{LP-F-G} (with $g=F$ and $p=\infty$, see \rf{G-LI1}) together with this inequality and \rf{N-WMP} implies the required inequality
 \begin{equation}\xlabel{INP-4}
\| F\| _{\LIR}\le C(m)\,\lambda
\end{equation}
 proving the lemma for $p=\infty$. Let us prove the lemma for  $p\in(1,\infty)$.
\smsk
\smallskip

We recall that $\{a_i\}_{i=-\infty}^{+\infty}$ is the strictly increasing sequence defined in \rf{A-SQN}. Let
 \begin{equation}\xlabel{L-Z}
L=\{i\in\Z:~a_i\in E\}.
\end{equation}

We know that $F\vert_{\tE}=\tf$, $\tf\vert  _{\tE\setminus E}\equiv 0$ and $\tf\vert  _E=f$. Thus, $F(a_i)=f(a_i)$ for each $i\in L$, and  $F(a_i)=0$ for each $i\in\Z\setminus L$. We apply Lemma \reff{LP-F-G} (taking $g=F$) and obtain the following:
 \begin{equation}\xlabel{FLPR}
\| F\| _{\LPR}^p\le C(m)^p\,\left\{\| F\| _{\LMPR}^p+B\right\}.
\end{equation}
 Here
 \begin{equation}\xlabel{B-D1}
B=\smed_{i\in L}\,\vert  f(a_i)\vert  ^p~~~\text{provided}~~~L\ne\emp,~~~\text{and}~~~B=0~~~\text{otherwise}.
\end{equation}

Inequalities \rf{N-WMP} and \rf{FLPR} imply that
 \begin{equation}\xlabel{FP-B}
\| F\| _{\LPR}^p\le C(m)^p\,\{\lambda^p+B\}.
\end{equation}

Prove that $B\le C(m)^p\,\lambda^p$. Indeed, this is trivial whenever $L=\emp$ (because $B=0$.) Assume that $L\ne\emp$. Without loss of generality, we may also assume that $L$ is finite.
\smallskip

Lemma \reff{SM-E} (with $k=0$) tells us that
$\vert  f(a_{i})\vert  \le C(m)\,\lambda$ for every $i\in L$,
so that $B\le C(m)^p\,\lambda^p$ holds provided $\#L\le m$.

Suppose that $\#L\ge m+1$. For simplicity of notation, we may assume that in this case $L=\{0, \ldots ,n\}$ where $n=\# L\ge m+1$. Thus
\[
B=\smed_{i=0}^n\,\vert  f(a_{i})\vert  ^p
\]
 where $\{a_i\}_{i=0}^n$ is a strictly increasing sequence in $E$ such that $a_{i+1}-a_i\ge 3$ for all $i=0, \ldots ,n-1$. See \rf{A-MD}. Hence,
\[
B=\smed_{i=0}^n\,\min\{1,a_{i+m}-a_i\}\,
\vert  \Delta^0f[a_i]\vert  ^p.
\]
 (We recall our Convention \reff{AGREE-2} which for the case of the sequence $\{a_i\}_{i=0}^n$ states that $a_i=+\infty$ provided $i>n$.) Therefore, thanks to Assumption \reff{ASMP-FE}, $B\le\,\lambda^p$.
\smsk

Thus, inequality $B\le C(m)^p\,\lambda^p$ holds, which together with \rf{FP-B} implies the statement of the lemma.\bx

\bigskip
\par {\it Proof of the sufficiency part of Theorem \reff{W-VAR-IN}.}~
It is well known that for every $p\in(1,\infty]$ we have
 \begin{equation}\xlabel{EQ-NW}
\| F\| _{\WMPR}=\smed_{k=0}^m\, \| F^{(k)}\| _{\LPR}\le C(m)\,(\| F\| _{\LPR}+\| F\| _{\LMPR}).
\end{equation}
  See, e.g.,~\cite{M}, p. 26.

We have proved that every function $f$ on $E$ can be extended to a function $F:\R\to\R$ such that
\[
\| F\| _{\LMPR}+\| F\| _{\LPR}\le C(m)\,\lambda
\]
 provided Assumption \reff{ASMP-FE} holds. See \rf{N-WMP} and Lemma \reff{LP-F}. Hence, thanks to \rf{EQ-NW}, we conclude that $F\in \WMPR$ and $\| F\| _{\WMPR}\le C(m)\,\lambda$.

Thus, $f\in\WMPR\vert  _E$ (because $F\vert_E=f$). Furthermore,
from definitions \rf{N-WMPR} and \rf{LMBD}, we have
\[
\| f\| _{\WMPR\vert  _E}\le \| F\| _{\WMPR}\le C(m)\,\lambda =C(m)\,\NWMP(f:E)
\]
 proving the sufficiency part of Theorem \reff{W-VAR-IN}.\bx
 
The proof of Theorem \reff{W-VAR-IN} is complete.\bx


\SECT{5. $\WMPR$-functions on sequences of points.}{5}
\addtocontents{toc}{5. $\WMPR$-functions on sequences of points. \hfill\thepage\par \VST}

\indent\par In this section we prove Theorem \reff{W-TFIN}. Let $p\in(1,\infty)$ and let $E=\{x_i\}_{i=\ell_1}^{\ell_2}$ be a strictly increasing sequence of points in $\R$ where
$\ell_1,\ell_2\in\Z\cup\{\pm\infty\}$, $\ell_1\le\ell_2$. In this case we assume that the following version of
convention \rf{AGR} holds.
\begin{convention}
\lbl{AGREE-3} $x_i=+\infty$ whenever $i>\ell_2$.
\end{convention}
      \smsk

In particular, this notational convention leads us to a certain simplification of the formula \rf{S-WP} for the quantity $\TLNW(f:E)$ whenever $E$ is a finite subset of $\R$ consisting of at most $m$ points.

Indeed, let  $0\le n<m$ and let $E=\{x_0, \ldots ,x_n\}$ where $x_0<\cdots<x_n$. In this case
\[
\ME=\min\{\meh m,\#E-1\}=n
\]
 (see \rf{ME}), so that
\[
\TLNW(f:E)=
\,\left(\,\smed_{k=0}^{n}\,\,\smed_{i=0}^{n-k}
\min\left\{1,x_{i+m}-x_{i}\right\}
\,\left\lvert  \Delta^kf[x_i, \ldots ,x_{i+k}]\right\rvert  ^p
\right)^{\frac1p}.
\]

Clearly, $i+m>n$ for each $0\le i\le n$ so that, according to Convention \reff{AGREE-3}, $x_{i+m}=+\infty$.  Therefore $\min\{1,x_{i+m}-x_{i}\}=1$ for every $i=0, \ldots ,n$, proving that
 \begin{equation}\xlabel{WN-RPR}
\TLNW(f:E)\sim
\,\max\{\,\vert  \Delta^kf[x_i, \ldots ,x_{i+k}]\vert  : k=0, \ldots ,n, \,i=0, \ldots ,n-k\}
\end{equation}
 with constants in this equivalence depending only on $m$.
 \msk
\smallskip

We turn to the proof of the necessity part of Theorem \reff{W-TFIN}.

\textit{(Necessity.)} For $p=\infty$ the necessity part of Theorem \reff{W-TFIN} is immediate from the necessity part of Theorem \reff{W-VAR-IN} (proven for an arbitrary closed $E\subset\R$). Let us prove the necessity for $p\in(1,\infty)$.

First, consider the case of a set $E\subset\R$ with $\#E\le m$. Let $0\le n<m$, and let $E=\{x_0, \ldots ,x_n\}$, $x_0<\cdots<x_n$. Let $f$ be a function on $E$, and let $F\in\WMPR$ be a function on $\R$ such that $F\vert_E=f$. Fix $k\in\{0, \ldots ,n\}$ and $i\in\{0, \ldots ,n-k\}$. Let us apply Lemma \reff{NP-W} to the function $G=F$, $I=\R$, $q=p$ and $S=\{x_i, \ldots ,x_{i+k}\}$. (Note that $F^{(k+1)}\in L_p(\R)$ because $F\in\WMPR$ and $k\le n\le m-1$). This lemma tells us that
\[
\vert  \Delta^{k}f[x_i, \ldots ,x_{i+k}]\vert  ^p\le 2^p\,\intl_{\R} \left(\vert  F^{(k)}(y)\vert  ^p+\vert  F^{(k+1)}(y)\vert  ^p\right)\,dy\le 2^{p+1}\,\| F\| _{\WMPR}^p
\]
 proving that
\[
\vert  \Delta^{k}f[x_i, \ldots ,x_{i+k}]\vert  \le 4\,\| F\| _{\WMPR}~~
\text{for all}~~
k\in\{0, \ldots ,n\}~~\text{and}~~i\in\{0, \ldots ,n-k\}.
\]

These inequalities and equivalence \rf{WN-RPR} imply that $\TLNW(f:E)\le C(m)\,\| F\| _{\WMPR}$. Taking the infimum in this inequality over all function $F\in\WMPR$ such that $F\vert_E=f$, we obtain the required inequality
 \begin{equation}\xlabel{TL-SME}
\TLNW(f:E)\le C(m)\,\| f\| _{\WMPR\vert  _E}.
\end{equation}

This proves the necessity in the case under consideration.
\smallskip

We turn to the case of a sequence $E$ containing at least $m+1$ elements. In this case
\[
\ME=\min\{\meh m,\,\#E-1\}=m~~~~\text{(see \rf{ME}).}
\]

The necessity part of Theorem \reff{W-VAR-IN} tells us that for every function $f\in\WMPR\vert  _E$ the following inequality
 \begin{equation}\xlabel{NM-F}
\NWMP(f:E)\le C(m)\,\| f\| _{\WMPR\vert  _E}
\end{equation}
 holds. We recall that $\NWMP(f:E)$ is the quantity defined by \rf{N-WRS-IN}.

Comparing \rf{N-WRS-IN} with \rf{S-WP}, we conclude that $\TLNW(f:E)\le\NWMP(f:E)$. This inequality and inequality \rf{NM-F} imply \rf{TL-SME} completing the proof of the necessity part of Theorem \reff{W-TFIN}.\bx
\medskip

\textit{(Sufficiency.)} We proceed the proof of the sufficiency by cases.
\msk

{\sc Case 1.} $p\in(1,\infty)$ \textit{and $\#E>m$ or $p=\infty$ and $E$ is an arbitrary sequence.}
\smsk

We recall that $E=\{x_i\}_{i=\ell_1}^{\ell_2}$ is a strictly increasing sequence in $\R$. If $p\in(1,\infty)$ then $E$ contains at least $m+1$ element, so that in this case $\ell_1+m\le \ell_2$. Furthermore, in this case $\ME=\min\{\meh m,\#E-1\}=m$.

Let $f$ be a function on $E$ satisfying the hypothesis of Theorem \reff{W-TFIN}, i.e.,
\[
\tlm=\TLNW(f:E)<\infty.
\]
This inequality and definition \rf{S-WP} enable us to make the following
\begin{assumption}
\lbl{S-FN} The following inequality
\[
\smed_{k=0}^{m}\,\,\smed_{i=\ell_1}^{\ell_2-k}
\min\left\{1,x_{i+m}-x_{i}\right\}
\,\left\lvert  \Delta^kf[x_i, \ldots ,x_{i+k}]\right\rvert  ^p\le \tlm^p
\]
 holds provided $1<p<\infty$.

Let $p=\infty$. Then for every $k\in\{0, \ldots ,\ME\}$ and every $i=\ell_1, \ldots ,\ell_2-k$ we have
 \begin{equation}\xlabel{A-SQ2}
\vert  \Delta^kf[x_i, \ldots ,x_{i+k}]\vert  \le \tlm.
\end{equation}
\end{assumption}

Our aim is to show that
\[
f\in\WMPR\vert  _E~~~~\text{and}~~~~\| f\| _{\WMPR\vert  _E}\le C(m)\,\tlm.
\]

We prove these properties of $f$ using a slight modification of the proof of the sufficiency part of Theorem \reff{W-VAR-IN} given in Sections 4. More specifically, we will only implement minor changes into Lemma \reff{SM-E} and Lemma \reff{LM-EW} related to replacing in their proofs the constant $\lambda$ with the constant $\tlm$, Assumption \reff{ASMP-FE} with Assumption \reff{S-FN}, and using tuples of \textit{consecutive} elements of the sequence $E=\{x_i\}_{i=\ell_1}^{\ell_2}$ (rather than arbitrary finite subsequences of $E$ which we used in Lemmas \reff{SM-E} and \reff{LM-EW}.)
\smsk

We begin with an analogue of Lemma \reff{SM-E}.
\begin{lemma}
 \lbl{DD-L} Let $k\in\{0, \ldots ,\ME-1\}$, $\nu\in\Z$, $\ell_1\le\nu\le\ell_2-k$, and let $\{y_i\}_{i=0}^n=\{x_{\nu+i}\}_{i=0}^n$ where $n=\ell_2-\nu$. Then
 \begin{equation}\xlabel{DK-T1}
\vert  \Delta^kf[y_0, \ldots ,y_k]\vert  \le C(m)\,\tlm\,.
\end{equation}
\end{lemma}

{\it Proof.} For $p=\infty$ the lemma is immediate from \rf{A-SQ2}.

Let $p\in(1,\infty)$. (Recall that in this case $\ME=m$.) Let $S=\{y_0, \ldots ,y_n\}$. If $y_m-y_0\ge 2$, then
\begin{align*}
A&=\vert  \Delta^kf[y_0, \ldots ,y_k]\vert  ^p=\min\{1,y_m-y_0\}\,
\vert  \Delta^kf[y_0, \ldots ,y_k]\vert  ^p\\
&=\min\{1,x_{\nu+m}-x_{\nu}\}\,
\vert  \Delta^kf[x_{\nu}, \ldots ,x_{\nu+k}]\vert  ^p
\end{align*}
 so that, thanks to Assumption \reff{S-FN}, $A\le\tlm^p$.
\smallskip

Now let $y_m-y_0<2$, and let $\ell$ be a positive integer, $m\le \ell\le n$, such that $y_\ell-y_0\le 2$ but $y_{\ell+1}-y_0>2$. (We recall that, according to Convention \reff{AGREE-3}, $x_i=+\infty$ whenever $i>\ell_2$, so that $y_i=+\infty$ if $i>n$.)

In this case Lemma \reff{P-WQ} tells us that
\begin{eqnarray}
A&\le& C(m)^p\,
\smed_{j=0}^{m}\,\,\smed_{i=0}^{\ell-j}
\min\left\{1,y_{i+m}-y_{i}\right\}
\,\left\lvert  \Delta^jf[y_i, \ldots ,y_{i+j}]\right\rvert  ^p\nn\\
&=&
C(m)^p\,
\smed_{j=0}^{m}\,\,\smed_{i=0}^{\ell-j}
\min\left\{1,x_{\nu+i+m}-x_{\nu+i}\right\}
\,\left\lvert  \Delta^jf[x_{\nu+i}, \ldots ,x_{\nu+i+m}]\right\rvert  ^p.
\nn
\end{eqnarray}
 This inequality and Assumption \reff{S-FN} imply \rf{DK-T1} proving the lemma.\bx

We turn to the analogue of Lemma \reff{LM-EW}. We recall that the set $G$ defined by \rf{G-WE} consists of isolated points of $\R$. (Moreover, the distance between any two distinct points of $G$ is at least $1$.) Since $E=\{x_i\}_{i=\ell_1}^{\ell_2}$ is a strictly increasing sequence of points, the set $\tE=E\cup G$ (see \rf{EW-G}) can be represented as a certain bi-infinite strictly increasing sequence of points:
 \begin{equation}\xlabel{TE-SQ1}
\tE=\{t_i\}_{i=-\infty}^{+\infty}.
\end{equation}

We also recall that the function $\tf:\tE\to\R$ is defined by formula \rf{TF-D2}.
\begin{lemma}
 \lbl{LL-2} Let $\ell\in\Z$, $n\in\N$, $n\ge m$, and let $y_i=t_{i+\ell}$, $i=0, \ldots ,n$. Then the following inequality
 \begin{equation}\xlabel{FW-NR}
\smed_{i=0}^{n-m}\,
(y_{i+m}-y_{i})\left\lvert  \Delta^m\tf[y_i, \ldots ,y_{i+m}]\right\rvert  ^p
\le\,C(m)^p\,\tlm^p
\end{equation}
 holds provided $1<p<\infty$. If $p=\infty$, then for every $i\in\Z$ we have
 \begin{equation}\xlabel{DPI-3}
\vert  \Delta^m\tf[t_i, \ldots ,t_{i+m}]\vert  \le\,C(m)\,\tlm.
\end{equation}
\end{lemma}

{\it Proof.} Let $1<p<\infty$. Repeating the proof of Lemma \reff{LM-EW} for the sequence $E=\{x_i\}_{i=\ell_1}^{\ell_2}$ we show that this lemma holds under a weaker hypothesis for the function $f$ on $E$. More specifically, we assume that $f$ satisfies the condition of Assumption \reff{S-FN} rather than Assumption \reff{ASMP-FE} (as for the case of an arbitrary closed set $E\subset\R$). In other words, we prove that everywhere in the proof of Lemma \reff{LM-EW} the constant $\lambda$ can be replaced with the constant $\tlm$ provided $E$ is a sequence and Assumption \reff{S-FN} holds for $f$.

The validity of such a replacement relies on the following obvious observation:

\textit{Let $0\le k\le m$, $0\le i\le n-k$, and let $t_i, \ldots ,t_{i+k}\in\tE$ be $k+1$ consecutive elements of the sequence $\tE$, see \rf{TE-SQ1}. If\,  $t_i, \ldots ,t_{i+k}$ belong to the sequence $E=\{x_i\}_{i=\ell_1}^{\ell_2}$, then $t_i, \ldots ,t_{i+k}$ are $k+1$ consecutive elements of this sequence. In other words, there exists $\nu\in\Z$, $\ell_1\le \nu\le\ell_2-k$, such that $t_i=x_{\nu+i}$ for all $i=0, \ldots ,k$.}
\smsk

In particular, this observation enables us to replace the constant $\lambda$ with $\tlm$ in inequalities \rf{L1}, \rf{L-A2}, \rf{L2}, and all inequalities after \rf{L2} until the end of the proof of Lemma \reff{LM-EW}.
\smsk

After all these modifications and changes, we literally follow the proof of Lemma \reff{LM-EW}. This leads us to the required inequality \rf{FW-NR} proving the lemma for $p\in(1,\infty)$.

For $p=\infty$ we repeat the proof given in Lemma \reff{LM-EW}. See \rf{PI-R}.

The proof of the lemma is complete.\bx

Because for $p\in(1,\infty)$ the integer $\ell$ from the hypothesis of Lemma \reff{LL-2} is \textit{arbitrary}, and the right hand side of inequality \rf{FW-NR} does not depend on $\ell$, we conclude that
\[
\smed_{i=-\infty}^{+\infty}\,
(t_{i+m}-t_{i})\left\lvert  \Delta^m\tf[t_i, \ldots ,t_{i+m}]\right\rvert  ^p
\le\,C(m)^p\,\tlm^p.
\]

This inequality and \rf{DPI-3} tell us that the function $\tf:\tE\to\R$ satisfies the hypothesis of Theorem \reff{DEBOOR}. This theorem produces a function $\tF\in\LMPR$ with $\| \tF\| _{\LMPR}\le C(m)\,\tlm$ such that $\tF\vert  _{\tE}=\tf$. Since $\tf\vert  _E=f$, the function $\tF$ is an extension of $f$ from $E$ to all of $\R$, i.e.,~$\tF\vert  _E=f$.

The following lemma is an analogue of Lemma \reff{LP-F} for the case of sequences.
\begin{lemma}
 $\| \tF\| _{\LPR}\le C(m)\,\tlm$\,\,  for every
$p\in(1,\infty]$.
\end{lemma}

{\it Proof.} For $p=\infty$ the proof literally follows the proof of the same case given in Lemma \reff{LP-F} (with corresponding replacing of $\lambda$ with $\tlm$). See \rf{INP-4}.
\smsk

We turn to the case $p\in(1,\infty)$ and $\#E\ge m+1$. For this case our proof relies on a slight modification of the proof given in Lemma \reff{LP-F}.

Let $\{a_i\}_{i=-\infty}^{+\infty}\subset\tE$ be the  bi-infinite strongly increasing sequence determined by \rf{A-SQN}, and let $B$ be the quantity defined by formula \rf{B-D1}. We also recall the definition of the family of indexes $L=\{i\in\Z:~a_i\in E\}$ given in \rf{L-Z}. Let $\vkp=\#L-1$.

We follow the proof of Lemma \reff{LP-F} and obtain an analogue of inequality \rf{FP-B} which states that
 \begin{equation}\xlabel{WFP-B}
\| \tF\| _{\LPR}^p\le C(m)^p\,\{\tlm^p+B\}
\end{equation}
   where $B$ is defined by \rf{B-D1}. Thus, our task is to prove that
 \begin{equation}\xlabel{B-A}
B=\smed_{i\in L}\,\vert  f(a_i)\vert  ^p\le C(m)^p\,\tlm^p.
\end{equation}

As in Lemma \reff{LP-F}, we may assume that $L\ne\emp$ (otherwise $B=0$). Since $\{a_i\}_{i=-\infty}^{+\infty}$ is strictly increasing, the family of points $\{a_i:i\in L\}$  is a subsequence of this sequence lying in the strictly increasing sequence $E=\{x_i\}_{i=\ell_1}^{\ell_2}$. This enables us to consider the family $\{a_i:i\in L\}$ as a strictly increasing subsequence of $\{x_i\}_{i=\ell_1}^{\ell_2}$, i.e.,~$a_i=x_{i_\nu}$ where $\nu=0, \ldots ,\vkp$, and
$i_{\nu_1}<i_{\nu_2}$ for all $0\le \nu_1<\nu_2\le \vkp$.
\smsk

Note that, according to Convention \reff{AGREE-3}, $x_{i_\nu+m}=+\infty$ provided $0\le\nu\le \vkp$ and  $i_\nu+m>\ell_2$. In particular, in this case $\min\{1,x_{i_\nu+m}-x_{i_\nu}\}=1$.
 \smsk

Let us partition the family $L$ into the following  two subfamilies:
\[
L_1=\{i\in L:~a_i=x_{i_\nu},~x_{i_{\nu}+m}-x_{i_\nu}\ge 2\},~
L_2=\{i\in L:~a_i=x_{i_\nu},~x_{i_{\nu}+m}-x_{i_\nu}<2\}.
\]

Clearly, for every $i\in L_1$,
\[
\vert  f(a_i)\vert  ^p=\vert  f(x_{i_\nu})\vert  ^p=
\min\{1,x_{i_\nu+m}-x_{i_\nu}\}\,\vert  f(x_{i_\nu})\vert  ^p.
\]
   Hence,
 \begin{equation}\xlabel{B1-Y}
B_1=\smed_{i\in L_1}\,\vert  f(a_i)\vert  ^p=
\smed_{\nu:\,i_\nu\in L_1}\,
\min\left\{1,x_{i_\nu+m}-x_{i_\nu}\right\}
\,\left\lvert  \Delta^0f[x_{i_\nu}]\right\rvert  ^p.
\end{equation}
 This inequality and Assumption \reff{S-FN} imply that
\[
B_1\le \,\smed_{i=\ell_1}^{\ell_2}
\min\left\{1,x_{i+m}-x_{i}\right\}
\,\left\lvert  \Delta^0f[x_i]\right\rvert  ^p\le \tlm^p.
\]

Let us estimate the quantity
 \begin{equation}\xlabel{B2-Y}
B_2=\smed_{i\in L_2}\,\vert  f(a_i)\vert  ^p.
\end{equation}

Lemma \reff{DD-L} tells us that
 $
\vert  f(a_i)\vert  =\vert  \Delta^0f[a_i]\vert  \le C(m)\,\tlm$ for every $i\in L_2$.
 Hence,
\[
B_2\le C(m)^p\,(\# L_2)\,\tlm^p\le C(m)^p\,(\# L)\,\tlm^p,
\]
 so that $B_2\le C(m)^p\,\tlm^p$ provided $\#L=\vkp+1\le m+1$.
\smsk
\smallskip

Prove that $B_2$ satisfies the same inequality whenever $\#L=\vkp+1>m+1$, i.e.,~$\vkp>m$. In particular, in this case \[\#E=\ell_2-\ell_1+1\ge \#L>m+1\] so that
$\ell_2-\ell_1>m$.

Fix $i\in L_2$. Thus $a_i=x_{i_\nu}$ for some $\nu\in\{0, \ldots ,\vkp\}$, and $x_{i_{\nu}+m}-x_{i_\nu}<2$. We know that
\[
x_{i_{\nu+1}}-x_{i_\nu}=a_{\nu+1}-a_{\nu}>2.
\]
 Therefore, there exists $\ell_\nu\in\N$, $i_{\nu}+m\le \ell_\nu< i_{\nu+1}$, such that $x_{\ell_{\nu}}-x_{i_\nu}\le 2$ but $x_{\ell_{\nu}+1}-x_{i_\nu}>2$.

Let $n=\ell_2-i_\nu$, and let $y_s=x_{i_\nu+s}$, $s=0, \ldots ,n$. We note that $n\ge m$; otherwise,  $i_{\nu}+m>\ell_2$ so that, according to Convention \reff{AGREE-3}, $x_{i_{\nu}+m}=+\infty$ which contradicts inequality $x_{i_{\nu}+m}-x_{i_\nu}<2$.

Let us apply Lemma \reff{P-WQ} to the strictly increasing sequence $\{y_s\}_{s=0}^n$ and a function $g(y_s)=f(x_{i_\nu+s})$ defined on the set $S=\{y_0, \ldots ,y_n\}$ with parameters $k=0$ and  $\ell=\ell_\nu$.

This lemma tells us that
\begin{eqnarray}
\vert  f(y_0)\vert  ^p=\vert  \Delta^0g[y_0]\vert  ^p&\le&  C(m)^p\,
\smed_{j=0}^{m}\,\,\smed_{s=0}^{\ell-j}
\min\left\{1,y_{s+m}-y_{s}\right\}
\,\left\lvert  \Delta^jg[y_s, \ldots ,y_{s+j}]\right\rvert  ^p\nn\\
&=&
C(m)^p\,
\smed_{j=0}^{m}\,\,\smed_{s=0}^{\ell_\nu-j}
\min\left\{1,x_{i_\nu+s+m}-x_{i_\nu+s}\right\}
\left\lvert  \Delta^jf[x_{i_\nu+s}, \ldots ,x_{i_\nu+s+j}]\right\rvert  ^p.
\nn
\end{eqnarray}

Because $\ell_\nu\le i_{\nu+1}-1$, we have
\begin{align*}
\vert  f(a_i)\vert  ^p&=\vert  f(y_0)\vert  ^p\\
&\le  C(m)^p\,
\smed_{j=0}^{m}\,\,\smed_{s=0}^{i_{\nu+1}-1-j}
\min\left\{1,x_{i_\nu+s+m}-x_{i_\nu+s}\right\}
\,\left\lvert  \Delta^jf[x_{i_\nu+s}, \ldots ,x_{i_\nu+s+j}]\right\rvert  ^p.
\end{align*}

Summarizing these inequalities over all $i\in L_2$, we conclude that
\begin{eqnarray}
B_2=\smed_{i\in L_2}\,\vert  f(a_i)\vert ^p
&\le&
C(m)^p\,\smed_{i=i_\nu\in L_2}\,
\smed_{j=0}^{m}\,\,\,\smed_{s=0}^{i_{\nu+1}-1-j}
\min\left\{1,x_{i_\nu+s+m}-x_{i_\nu+s}\right\}\nn\\
&\times&\,\left\lvert  \Delta^jf[x_{i_\nu+s}, \ldots ,x_{i_\nu+s+j}]\right\rvert  ^p
\nn\\
&=&
C(m)^p\,
\smed_{j=0}^{m}\,\,\,
\smed_{i=i_\nu\in L_2}\,\,
\smed_{s=0}^{i_{\nu+1}-1-j}
\min\left\{1,x_{i_\nu+s+m}-x_{i_\nu+s}\right\}\nn\\
&\times&\,\left\lvert  \Delta^jf[x_{i_\nu+s}, \ldots ,x_{i_\nu+s+j}]\right\rvert  ^p
\nn\\
&\le&
C(m)^p\,
\smed_{j=0}^{m}\,\,\smed_{i=\ell_1}^{\ell_2-j}
\min\left\{1,x_{i+m}-x_{i}\right\}
\,\left\lvert  \Delta^jf[x_i, \ldots ,x_{i+j}]\right\rvert  ^p,
\nn
\end{eqnarray}
 so that, thanks to Assumption \reff{S-FN}, $B_2\le C(m)^p\,\tlm^p$.
 \msk

\smallskip
Finally, from \rf{B-A}, \rf{B1-Y} and \rf{B2-Y}, we have
\[
B=B_1+B_2\le \tlm^p+C(m)^p\,\tlm^p
\]
 proving \rf{B-A}. This inequality together with \rf{WFP-B} completes the proof of the lemma.\bx
       \msk

Thus, $\| \tF\| _{\LMPR}\le C(m)\,\tlm$ and $\| \tF\| _{\LPR}\le C(m)\,\tlm$. These inequalities together with \rf{EQ-NW} imply that
\[
\| \tF\| _{\WMPR}\le C(m)\,(\| \tF\| _{\LPR}+\| \tF\| _{\LMPR})\le
C(m)\,\tlm.
\]

Because $\tF\in\WMPR$ and $\tF\vert  _E=f$, the function $f$ belongs to the trace space $\WMPR\vert  _E$. Furthermore,
\[
\| f\| _{\WMPR\vert  _E}\le \| \tF\| _{\WMPR}\le C(m)\,\tlm =C(m)\,\TLNW(f:E)
\]
 proving the sufficiency part of Theorem \reff{W-TFIN} in the case under consideration.
\bsk

\smallskip
{\sc Case 2.} $p\in(1,\infty)$ \textit{and $\#E\le m$.}
\smsk

We may assume that $E=\{x_i\}_{i=0}^n$ is a strictly increasing sequence, and $0\le n<m$. (In other words, we assume that $\ell_1=0$ and $\ell_2=n$.) In this case $\ME=\min\{\meh m,\,\#E-1\}=n$ (see \rf{ME}), and equivalence \rf{WN-RPR} holds with constants depending only on $m$.

Let $f$ be a function on $E$, and let $\tlm=\TLNW(f:E)$. In these settings, equivalence \rf{WN-RPR} enables us to make the following
\begin{assumption}
\lbl{SM-FNQ} For every $k=0, \ldots ,n$, and every $i=0, \ldots ,n-k$, the following inequality
\[
\vert  \Delta^kf[x_i, \ldots ,x_{i+k}]\vert  \le C(m)\,\tlm
\]
 holds.
\end{assumption}

Our aim is to prove the existence of a function $\tF\in\WMPR$ with $\| \tF\| \le C(m)\,\tlm$ such that $\tF\vert  _E=f$. We will do this by reduction of the problem to the {\sc Case 1}. More specifically, we introduce $m-n$ additional points $x_{n+1},x_{n+2}, \ldots , x_m$ defined by
 \begin{equation}\xlabel{NP-E}
x_k=x_n+2(k-n), ~~~~k=n+1, \ldots ,m.
\end{equation}

Let $\Ec=\{x_0, \ldots ,x_n,x_{n+1}, \ldots ,x_m\}$. Clearly, $\Ec$ is a strictly increasing sequence in $\R$ with $\#\Ec=m+1$. Furthermore, $E\subset \Ec$ and $E\ne \Ec$.

By $\brf:\Ec\to \R$ we denote \textit{the extension of $f$ from $E$ to $\Ec$ by zero}. Thus,
 \begin{equation}\xlabel{BRF-DF}
\brf(x_i)=f(x_i)~~\text{for all}~~0\le i\le n, ~\text{and}~~\brf(x_i)=0~~\text{for all}~~n+1\le i\le m.
\end{equation}

Then we prove that the following analogue of Assumption \reff{S-FN} (with $\ell_1=0$ and $\ell_2=m$) holds for the sequence $\Ec=\{x_i\}_{i=0}^m$:
 \begin{equation}\xlabel{A-7}
\TLNW(\brf:\Ec)^p=\smed_{k=0}^{m}\,\,\smed_{i=0}^{m-k}
\min\left\{1,x_{i+m}-x_{i}\right\}
\,\left\lvert  \Delta^k\brf[x_i, \ldots ,x_{i+k}]\right\rvert  ^p\le C(m)^p\,\tlm^p.
\end{equation}

Let us note that, according to Convention \reff{AGREE-3}, $x_{i+m}=+\infty$ provided $i=1, \ldots ,n$. Furthermore, thanks to \rf{NP-E}, $x_m-x_0=x_n+2(m-n)-x_0\ge 2$. Hence,
 \begin{equation}\xlabel{W-ND}
\TLNW(\brf:\Ec)^p=\smed_{k=0}^{m}\,\,\smed_{i=0}^{m-k}
\,\left\lvert  \Delta^k\brf[x_i, \ldots ,x_{i+k}]\right\rvert  ^p.
\end{equation}

Let us prove that
 \begin{equation}\xlabel{DD-8}
\vert  \Delta^k\brf[x_i, \ldots ,x_{i+k}]\vert  \le C(m)\,\tlm
\end{equation}
 for every $k=0, \ldots ,m$ and every $i=0, \ldots ,m-k$.
\smsk

\smallskip
Let $S=\{x_i, \ldots ,x_{i+k}\}$. Assumption \reff{SM-FNQ} tells us that \rf{DD-8} holds for every $k=0, \ldots ,n$ and every $i=0, \ldots ,n-k$. It is also clear that \rf{DD-8} holds for every $k, 0\le k<m-n$, and every $i, n<i\le m-k$, because in this case $\brf\vert  _S\equiv 0$. See \rf{BRF-DF}.

Let $0\le i\le n$ and let $n<i+k\le m$. Then, thanks to \rf{NP-E}, $x_{i+k}-x_i\ge x_{n+1}-x_n=2$ proving that $\diam S\ge 2$. This enables us to apply Lemma \reff{SM-DF} to  $S$ and the restriction of the function $\brf$ to $S$. Lemma \reff{SM-DF} tells us that there exists
$j\in\{0, \ldots ,k-1\}$ and $\nu\in\{i, \ldots ,i+k-j\}$ such that $x_{\nu+j}-x_\nu\le 1$ and
\[
\vert  \Delta^k\brf[x_i, \ldots ,x_{i+k}]\vert  \le\,2^k \,\vert  \Delta^j\brf[x_\nu, \ldots ,x_{\nu+j}]\vert  /\diam S.
\]

Because $\diam S\ge 2$ and $0\le k\le m$, we have
 \begin{equation}\xlabel{NU-LE}
\vert  \Delta^k\brf[x_i, \ldots ,x_{i+k}]\vert  \le\,2^m \,\vert  \Delta^j\brf[x_\nu, \ldots ,x_{\nu+j}]\vert  \,.
\end{equation}

Prove \rf{DD-8} for the case $j=0$. In this case
 \begin{equation}\xlabel{JZA}
\vert  \Delta^j\brf[x_\nu, \ldots ,x_{\nu+j}]\vert  =
\vert  \Delta^0\brf[x_\nu]\vert  =\vert  \brf(x_\nu)\vert  .
\end{equation}
 We know that $\brf(x_\nu)=0$ for $n+1\le \nu\le m$, see
\rf{BRF-DF}. In turn, Assumption \reff{SM-FNQ} (with $k=0$) and definition \rf{BRF-DF} tell us that
$\vert  \brf(x_\nu)\vert  =\vert  f(x_\nu)\vert  \le C(m)\,\tlm$ provided $0\le\nu\le n$.

From this inequality, \rf{NU-LE} and \rf{JZA}, we have
\[
\vert  \Delta^k\brf[x_i, \ldots ,x_{i+k}]\vert  \le
\,2^m\,\vert  \brf(x_\nu)\vert  \le\,2^m C(m)\,\tlm
\]
 proving \rf{DD-8} for $j=0$.
\smsk

Now let $0<j\le k-1$. In this case, from \rf{NP-E} and the inequality $x_{\nu+j}-x_\nu\le 1$ it follows that  $0<j\le n$ and $0\le \nu\le n-j$ (otherwise $x_{\nu+j}-x_\nu\ge 2$). Hence, thanks to \rf{BRF-DF},
$\brf(x_\ell)=f(x_\ell)$ for each $\ell=\nu, \ldots ,\nu+j$, so that
\[
\Delta^j\brf[x_\nu, \ldots ,x_{\nu+j}]=
\Delta^jf[x_\nu, \ldots ,x_{\nu+j}].
\]

This property and Assumption \reff{SM-FNQ} imply that
\[
\vert  \Delta^j\brf[x_\nu, \ldots ,x_{\nu+j}]\vert  =
\vert  \Delta^jf[x_\nu, \ldots ,x_{\nu+j}]\vert  \le C(m)\,\tlm
\]
 which together with \rf{NU-LE} yields the required inequality \rf{DD-8} in the case under consideration.
\smsk

\smallskip
We have proved that inequality \rf{DD-8} holds
for all $k=0, \ldots ,m$ and all $i=0, \ldots ,m-k$. The required inequality \rf{A-7} is immediate from \rf{DD-8} and formula \rf{W-ND}.
\smsk

Thus, $\#\Ec=m+1>m$ and an analogue of Assumption \reff{S-FN}, inequality \rf{A-7}, holds for the function $\brf$ defined on $\Ec$. This enables us to apply to $\Ec$ and $\brf$ the result of {\sc Case 1} which produces a function $\tF\in\WMPR$ such that $\tF\vert  _{\Ec}=\brf$ and
$\| \tF\| _{\WMPR}\le C(m)\,\tlm$.

Clearly, $\tF\vert  _{E}=f$ (because $\brf\vert  _E=f$), so that $f\in\WMPR\vert  _E$ and $\| f\| _{\WMPR\vert  _E}\le C(m)\,\tlm$ proving the sufficiency part of Theorem \reff{W-TFIN} in {\sc Case 2.}
\smsk

\smallskip
The proof of Theorem \reff{W-TFIN} is complete.\bx


\SECT{6. $W^m_p$-extension criteria in terms of local sharp maximal functions.}{6}
\addtocontents{toc}{6. $W^m_p$-extension criteria in terms of local sharp maximal functions. \hfill\thepage\par \VST}

\indent\par In this section we prove Theorem \reff{W-MF}.
\msk

\textit{(Necessity.)} Let $p\in(1,\infty)$, and let $F\in\WMPR$ be a function such that $F\vert_E=f$. Prove that for every $x\in\R$ and every $k=0, \ldots ,m$, the following inequality
 \begin{equation}\xlabel{FK-MF}
\fks(x)\le 4\,\smed_{j=0}^m\,\Mc[F^{(j)}](x)
\end{equation}
 holds. (Recall that $\Mc[g]$ denotes the Hardy--Littlewood maximal function. See \rf{HL-M}.)
\smsk

Let $0\le k\le m-1$, and let $S=\{x_0, \ldots ,x_k\}$, $x_0<\cdots<x_k$, be a $(k+1)$-point subset of $E$ such that $\diam(S\cup\{x\})\le 1$. Thus, $S\subset I=[x-1,x+1]$.

Let us apply Lemma \reff{NP-W} to the function $G=F$, the interval $I=[x-1,x+1]$ (with $\vert  I\vert  =2$), $q=1$ and the set $S$. (We note that $F^{(k+1)}\in L_p(I)$ because $F\in\WMPR$ and $k\le m-1$). This lemma tells us that
\[
\vert  \Delta^kF[S]\vert  \le 2\,\intl_I \left(\vert  F^{(k)}(y)\vert  +\vert  F^{(k+1)}(y)\vert  \right)\,dy.
\]
 Hence,
\begin{align*}
\vert  \Delta^kF[S]\vert  &\le 4\,\frac{1}{\vert  I\vert  }\,\intl_I \vert  F^{(k)}(y)\vert  \,dy+4\,\frac{1}{\vert  I\vert  }\,\intl_I \vert  F^{(k+1)}(y)\vert  \,dy\\
&\le 4\,\{\Mc[F^{(k)}](x)+\Mc[F^{(k+1)}](x)\}.
\end{align*}
 This inequality and definition \rf{FK-1} imply \rf{FK-MF} in the case under consideration.

Let us prove \rf{FK-MF} for $k=m$. Let $x\in\R$ and let $S=\{x_0, \ldots ,x_m\}$, $x_0<\cdots<x_m$, be an $(m+1)$-point subset of $E$ such that $\diam(S\cup\{x\})\,\le 1$. Let $I$ be the smallest closed interval containing $S\cup\{x\}$. Clearly, $\vert  I\vert  =\diam (S\cup\{x\})$. This and inequality \rf{DVD-IN} imply that
\begin{eqnarray}
\frac{\diam S}{\diam (S\cup\{x\})}\,\vert  \Delta^mf[S]\vert
&=&\frac{x_m-x_0}{\vert  I\vert  }\,\vert  \Delta^mf[S]\vert  \nn\\
&\le& \frac{x_m-x_0}{\vert  I\vert  }\cdot
\frac{1}{x_m-x_0}
\,\intl_{x_0}^{x_m}\,\vert  F^{(m)}(t)\vert  \,dt
\le \frac{1}{\vert  I\vert  }\,\intl_I\,\vert  F^{(m)}(t)\vert  \,dt.
\nn
\end{eqnarray}
 From this inequality and definition \rf{FK-2} it follows  that $\fms(x)\le \Mc[F^{(m)}](x)$ proving \rf{FK-MF} for $k=m$.

Now inequality \rf{FK-MF} and the Hardy--Littlewood maximal theorem yield
\begin{align*}
\smed_{k=0}^m\,\| \fks\| _{\LPR}&\le 4(m+1)\, \smed_{k=0}^m\,\| \Mc[F^{(k)}]\| _{\LPR}\\
&\le
C(m,p)\,\smed_{k=0}^m\,\,\| F^{(k)}\| _{\LPR}=
C(m,p)\,\| F\| _{\WMPR}.
\end{align*}
 Taking the infimum in the right hand side of this inequality over all $F\in\WMPR$ such that $F\vert_E=f$, we obtain the required inequality
\[
\smed_{k=0}^m\,\| \fks\| _{\LPR}\le
C(m,p)\,\| f\| _{\WMPR\vert  _E}
\]
 proving the necessity part of Theorem \reff{W-MF}.\bx
\msk
\smallskip

\textit{(Sufficiency.)} Let $f$ be a function on $E$ such that $\fks\in\LPR$ for every $k=0, \ldots ,m$, and let
\[
\WCP(f:E)=\smed_{k=0}^m\,\| \fks\| _{\LPR}.
\]
 Let us prove that $f\in\WMPR\vert  _E$ and $\| f\| _{\WMPR\vert  _E}\le C(m)\,\WCP(f:E)$.
\smsk

\smallskip
We begin with the case of a set $E$ with $\#E\ge m+1$. Let us show that in this case
 \begin{equation}\xlabel{W-4}
\NWMP(f:E)\le C(m)\,\WCP(f:E)
\end{equation}
 where $\NWMP$ is the quantity defined by \rf{N-WRS-IN}.
\smsk

Let $n\in\N$, $n\ge m$, and let $S=\{x_0, \ldots ,x_n\}$, $x_0<\cdots<x_n$, be a subset of $E$. Fix $k\in\{0, \ldots ,m\}$ and prove that for every $i\in\{0, \ldots ,n-k\}$ the following inequality
 \begin{equation}\xlabel{QH-1}
\min\left\{1,x_{i+m}-x_{i}\right\}
\,\left\lvert  \Delta^kf[x_i, \ldots ,x_{i+k}]\right\rvert  ^p\le
2^{mp}\, \smed_{j=0}^m\,\intl_{x_i}^{x_{i+m}}\,(\fjs)^p(u)\,du
\end{equation}
 holds.
 \smsk

We proceed by cases. Let $S_i=\{x_i, \ldots ,x_{i+k}\}$, and let $t_i=\min\left\{1,x_{i+m}-x_{i}\right\}$.
 \smsk
\smallskip

\textit{Case A.} $\diam S_i=x_{i+k}-x_i\le t_i$.
\smsk

Let $V_i=[x_i,x_i+t_i]$. Then $\vert  V_i\vert  \le 1$ and $V_i\supset[x_i,x_{i+k}]\supset S_i$ so that
$\diam(S_i\cup\{x\})\le 1$ for every $x\in V_i$. Therefore, thanks to \rf{FK-1},
 \begin{equation}\xlabel{KA-1}
\vert  \Delta^kf[S_i]\vert  \le \fks(x)~~~\text{for every}~~~x\in V_i~~~\text{and every}~~~k=0, \ldots ,m-1.
\end{equation}

In turn, if $k=m$, then
$S_i=\{x_i, \ldots ,x_{i+m}\}\subset V_i\subset [x_i,x_{i+m}]$ so that
\[
\diam S_i=\diam(S_i\cup\{x\})=x_{i+m}-x_i~~~~\text{for every}~~~~x\in V_i.
\]
 This and \rf{FK-2} imply that
\[
\vert  \Delta^mf[S_i]\vert  \le
\,\frac{\diam (S_i\cup\{x\})}{\diam S_i}\cdot \fms(x)=\fms(x),~~~~x\in V_i,
\]
 proving that inequality \rf{KA-1} holds for all $k=0, \ldots ,m$.

Raising both sides of \rf{KA-1} to the power $p$ and then integrating on $V_i$ with respect to $x$, we obtain that for each $k=0, \ldots ,m$, the following inequality
\[
\vert  V_i\vert
\left\lvert  \Delta^kf[x_i, \ldots ,x_{i+k}]\right\rvert  ^p
=\min\left\{1,x_{i+m}-x_{i}\right\}
\left\lvert  \Delta^kf[x_i, \ldots ,x_{i+k}]\right\rvert  ^p
\le
\intl_{x_i}^{x_{i+m}}(\fks)^pdx
\]
holds. Of course, this inequality implies \rf{QH-1} in the case under consideration.

\smallskip
\textit{Case B.} $\diam S_i=x_{i+k}-x_i> t_i$.
 \smsk

Clearly, in this case $t_i=1$ (because $x_{i+k}-x_i\le x_{i+m}-x_i$) so that $\diam S_i=x_{i+k}-x_i>1$.

In particular, $k\ge 1$. Lemma \reff{SM-DF} tells us that there exist $j\in\{0, \ldots ,k-1\}$ and $\nu\in\{i, \ldots ,i+k-j\}$ such that
$x_{\nu+j}-x_\nu\le 1$ and
 \begin{equation}\xlabel{DK-J}
\vert  \Delta^kf[S_i]\vert  \le\,2^k \,\vert  \Delta^jf[x_\nu, \ldots ,x_{\nu+j}]\vert  /\diam S_i
\le\,2^m \,\vert  \Delta^jf[x_\nu, \ldots ,x_{\nu+j}]\vert  \,.
\end{equation}

Because $x_{\nu+j}-x_\nu\le 1$,\, $x_{i+k}-x_i>1$ and $[x_\nu,x_{\nu+j}]\subset[x_i,x_{i+k}]$, there exists an interval $H_i$ such that $\vert  H_i\vert  =1$ and $[x_\nu,x_{\nu+j}]\subset H_i\subset[x_i,x_{i+k}]$.

Clearly, the set $\tS=\{x_\nu, \ldots ,x_{\nu+j}\}\subset H_i$ so that $\diam(\tS\cup\{x\})\le 1$ for every $x\in H_i$ (because $\vert  H_i\vert  =1$). Therefore, thanks to \rf{FK-1},~
$\vert  \Delta^jf[\tS]\vert  \le \fjs(x)$ for every $x\in H_i$, proving that
\[
\vert  \Delta^kf[S_i]\vert  \le\,2^m\,\fjs(x)~~~~\text{for all}~~~~x\in H_i.
\]
 See \rf{DK-J}. Raising both sides of this inequality to the power $p$ and then integrating on $H_i$ with respect to $x$, we obtain the following:
\begin{eqnarray}
\min\left\{1,x_{i+m}-x_{i}\right\}
\,\left\lvert  \Delta^kf[x_i, \ldots ,x_{i+k}]\right\rvert  ^p
&=&
\left\lvert  \Delta^kf[x_i, \ldots ,x_{i+k}]\right\rvert  ^p
=\vert  H_i\vert  \, \left\lvert  \Delta^kf[S_i]\right\rvert  ^p\nn\\
&\le&
2^{mp}\,\intl_{H_i}\,(\fjs)^p(x)\,dx\nn\\
&\le& 2^{mp}
\intl_{x_i}^{x_{i+m}}\,(\fjs)^p(x)\,dx.
\nn
\end{eqnarray}

This inequality implies \rf{QH-1} in \textit{Case B}
proving that \rf{QH-1} holds in both case.
\msk
\smallskip

Inequality \rf{QH-1} enables us to show that for every $k=0, \ldots ,m$, we have
 \begin{equation}\xlabel{AK-3}
A_k(f:S)=\smed_{i=0}^{n-k}
\min\left\{1,x_{i+m}-x_{i}\right\}
\,\left\lvert  \Delta^kf[x_i, \ldots ,x_{i+k}]\right\rvert  ^p\le C(m)^p\,\WCP(f:E)^p.
\end{equation}

 To prove \rf{AK-3}, we introduce a family of closed intervals
\[
\Tc=\{T_i=[x_i,x_{i+m}]:i=0, \ldots ,n-k\}.
\]

Note that each interval $T_{i_0}=[x_{i_0},x_{{i_0}+m}]\in\Tc$ has common points with at most $2m+1$ intervals from the family $\Tc$. In this case Lemma \reff{GRAPH} tells us that there exist a positive integer $\varkappa\le 2m+2$ and subfamilies $\Tc_\ell\subset \Tc$, $\ell=1, \ldots ,\varkappa$, each consisting of pairwise disjoint intervals, such that $\Tc=\cup\{\Tc_\ell:\ell=1, \ldots ,\varkappa\}$.

This property and \rf{QH-1} imply that
\[
A_k(f:S)\le 2^{mp}\,\smed_{\ell=1}^\varkappa\,
A_{k,\ell}(f:S)
\]
 where
\[
A_{k,\ell}(f:S)=\smed_{i:T_i\in\,\Tc_\ell}\,
\smed_{j=0}^m\,
\intl_{T_i}\,(\fjs)^p(u)\,du.
\]

Let $1\le \ell\le\varkappa$, and let $U_\ell=\cup\{T_i:\,T_i\in\Tc_\ell\}$. Because the intervals of each subfamily $\Tc_\ell$ are pairwise disjoint,
\[
A_{k,\ell}(f:S)=\smed_{j=0}^m\,
\intl_{U_\ell}\,(\fjs)^p(u)\,du\le
\smed_{j=0}^m\,
\intl_{\R}\,(\fjs)^p(u)\,du=
\smed_{j=0}^m\,\| \fjs\| ^p_{\LPR}.
\]
 Hence,
\[
A_k(f:S)\le 2^{mp}\,\smed_{\ell=1}^\varkappa\,
A_{k,\ell}(f:S)
\le 2^{mp}\,\varkappa\,\smed_{j=0}^m\,\| \fjs\| ^p_{\LPR}
\le (2m+2)\,2^{mp}\,\,\WCP(f:E)^p
\]
 proving \rf{AK-3}. From this inequality it follows that
\[
\smed_{k=0}^{m}\,\,\smed_{i=0}^{n-k}
\min\left\{1,x_{i+m}-x_{i}\right\}
\left\lvert  \Delta^kf[x_i, \ldots ,x_{i+k}]\right\rvert  ^p
= \smed_{k=0}^{m}\,A_k(f:S)\le
C^p\,\WCP(f:E)^p
\]
with $C=C(m)$. Taking the supremum in the left hand side of this inequality over all $n\in\N, n\ge m$, and all strictly increasing sequences $S=\{x_0, \ldots ,x_n\}\subset E$, and recalling definition \rf{N-WRS-IN}, we obtain the required inequality \rf{W-4}.
 \smsk
\smallskip

We know that $\fks\in\LPR$ for every $k=0, \ldots ,m$, so that $\WCP(f:E)<\infty$. This and inequality \rf{W-4} imply that $\NWMP(f:E)<\infty$ proving that $f$ satisfies the hypothesis of Theorem \reff{W-VAR-IN}. This theorem tells us that $f\in\WMPR\vert  _E$ and $\| f\| _{\WMPR\vert  _E}\le C(m)\,\NWMP(f:E)$.

\smallskip
This inequality together with \rf{W-4}
implies that
$
\| f\| _{\WMPR\vert  _E}\le C(m)\,\WCP(f:E)
$
 proving the sufficiency for every set $E\subset\R$ containing at least $m+1$ points.
\smsk
\medskip

It remains to prove the sufficiency for an arbitrary set $E\subset\R$ containing at most $m$ points.

 Let $0\le n<m$, and let $E=\{x_0, \ldots ,x_n\}$, $x_0<\cdots<x_n$, so that $\ME=\min\{\meh m,\,\#E-1\}=n$, see \rf{ME}. In this case Theorem \reff{W-TFIN} tells us that
$$\| f\| _{\WMPR\vert  _E}\le C(m)\,\TLNW(f:E)$$ where
\[
\TLNW(f:E)=
\,\left(\,\smed_{k=0}^{n}\,\,\smed_{i=0}^{n-k}
\min\left\{1,x_{i+m}-x_{i}\right\}
\,\left\lvert  \Delta^kf[x_i, \ldots ,x_{i+k}]\right\rvert  ^p
\right)^{\frac1p}.
\]
 We recall that $x_{i+m}=+\infty$ for all $i=0, \ldots ,n$ (see Convention \reff{AGREE-3}), so that in the above formula for $\TLNW(f:E)$ we have $\min\{1,x_{i+m}-x_{i}\}=1$ for every $i=0, \ldots ,n$.

Then we literally follow the proof of \rf{W-4}, and obtain the following inequality:
\[
\TLNW(f:E)\le C(m)\,\WCP(f:E).
\]

Hence, $\| f\| _{\WMPR\vert  _E}\le C(m)\,\WCP(f:E)$ proving the sufficiency for every set $E\subset\R$ consisting of at most $m$ points.
\smsk
\smallskip

The proof of Theorem \reff{W-MF} is complete.\bx
\bigskip

\textit{Further remarks and comments.}
\smallskip

$(\bigstar 1)$~Let $\EXT_E(\cdot:\WMPR)$ be the extension operator for the Sobolev space $\WMPR$ constructed in Section 4. (See \rf{EXT-W}). Examining this extension construction, we conclude that for each $f$ on $E$ the support of its extension $F=\EXT_E(f:\WMPR)$ has the  following useful property:
\begin{statement}
 Let $f$ be a function on $E$ such that $\NWMP(f:E)<\infty$.
Then the support of the function  $F=\EXT_E(f:\WMPR)$ lies in a $\delta$-neighborhood of $E$ where $\delta=3(m+2)$.
\end{statement}

For a proof of this statement we refer the reader
to~\cite[Section 5.5]{Sh-LV-2018}.
\smallskip

Let $0<\ve\le 1$ and let $\gamma\ge 1$. We can slightly modify
the construction of the extension $F$ from formula
\rf{EXT-W} by replacing the constant $n_J$ from  the
formula \rf{NJ-12} with the constant $n_J=\lfloor\vert  J\vert  /(\gamma\ve)\rfloor$. Obvious
changes in the proof of Theorem \reff{W-VAR-IN} enable us to show
that after such a modification this theorem holds with constants
in equivalence \rf{NM-CLC-IN} depending only on $m$,
$\gamma$ and $\ve$. Following the same scheme as in
the proof of~\cite[Statement 5.20]{Sh-LV-2018}, one can show that
the support of this modified extension $F$ lies in the
$\ve$-neighborhood of $E$ provided
$\gamma=\gamma(m)$ is big enough.

We leave the details to the interested reader.
       \msk
\smallskip

$(\bigstar 2)$~Let $E=\{x_i\}_{i=\ell_1}^{\ell_2}$ be a strictly increasing
sequence of points in $\R$, and let $f$ be a
function on $E$. As we note in~\cite{Sh-LV-2018}, in
this case the extension method suggested in this paper, produces
a piecewise polynomial $C^{m-1}$-smooth extension
$F=\EXT_E(f:\LMPR)$ which coincides with a polynomial of degree at most
$2m-1$ on each interval $(x_i,x_{i+1})$. (Details are
spelled out in~\cite[Section 4.4]{Sh-LV-2018}.)
\smsk

We note that the same statement holds for the extension operator $\EXT_E(\cdot:\WMPR)$ introduced in the present paper. This is immediate from formula \rf{EXT-WL} and the following obvious observation: \textit{the set $\tE$ (see \rf{EW-G}) is a sequence of points in $\R$ provided the set $E$ is}.

This enables us to reformulate this property of the extension $F=\EXT_E(f:\WMPR)$ in terms of Spline Theory as follows: \textit{The extension $F$ is an interpolating $C^{m-1}$-smooth spline of order\, $2m$\, with knots\, $\{x_i\}_{i=\ell_1}^{\ell_2}$.}
\smallskip

$(\bigstar 3)$~Let $0\le n<m$ and let $E=\{x_0, \ldots ,x_n\}$ where $x_0<\cdots<x_n$. Thus, in this case the constant $\ME=\min\{\meh m,\#E-1\}=n$, see \rf{ME}. Equivalence \rf{WN-RPR} and Theorem \reff{W-TFIN} tell us that for every $p\in(1,\infty)$ we have
 \begin{equation}\xlabel{SN-SQ}
\| f\| _{\WMPR\vert  _E}\sim \,\max\{\,\vert  \Delta^kf[x_i, \ldots ,x_{i+k}]\vert  : k=0, \ldots ,n, \,i=0, \ldots ,n-k\}
\end{equation}
 with constants in this equivalence depending only on $m$.
On the other hand, Theorem \reff{W-TFIN} also tells us that equivalence \rf{SN-SQ} holds for $p=\infty$ (again, with constants depending only on $m$).

Hence,
 \begin{equation}\xlabel{WRI-E}
\WMPR\vert  _E=W^n_\infty(\R)\vert  _E~~~\text{provided}~~~n=\#E-1<m,
\end{equation}
 with imbedding constants depending only on $m$.

In connection with isomorphism \rf{WRI-E}, we recall a classical Sobolev imbedding theorem which states that
 \begin{equation}\xlabel{AGV}
\WMPR\hookrightarrow  W^n_\infty(\R)~~~~\text{but}~~~~
\WMPR\ne W^n_\infty(\R).
\end{equation}

Comparing \rf{AGV} with \rf{WRI-E}, we observe that, even though the spaces $\WMPR$ and $W^n_\infty(\R)$ are distinct, their traces to every ``small'' subset of $\R$ containing at most $m$ points, coincide with each other.

\bigskip
\par\noindent {\bf Acknowledgements}
\bsk
\par I am very thankful to Michael Cwikel for useful suggestions and remarks. I am grateful to Charles Fefferman, Bo'az Klartag and Yuri Brudnyi for valuable conversations. I thank the referees for very careful reading and numerous comments and suggestions, which led to improvements of the manuscript.
\par The results of this paper were presented at the 11th Whitney Problems Workshop, Trinity College Dublin, Dublin, Ireland. I am very thankful to all participants of this conference for stimulating discussions and valuable advice.
\par This research was supported by Grant No 2014055 from the United States-Israel Binational Science Foundation (BSF).
\newpage


\begin{thebibliography}{ABCD}

\addtocontents{toc}{References \hfill \thepage\par}
{\footnotesize
\bibitem {deB4} C. de Boor, How small can one make the derivatives of an interpolating function? J. Approx.
Theory, 13 (1975) 105--116.
\bibitem {C1} A. P. Calder\'{o}n, Estimates for
    singular integral operators in terms of maximal
    functions, Studia Math. 44 (1972) 563--582.
\bibitem {CS} A. P. Calder\'{o}n, R. Scott, Sobolev
    type inequalities for $p>0$,  Studia Math. 62 (1978)
    75--92.
\bibitem {DL} R. DeVore, G. Lorentz, Constructive approximation, Fundamental Principles of Mathema\-ti\-cal Sciences, 303. Springer-Verlag, Berlin, 1993.
\bibitem {Es} D. Est\'evez, Explicit traces of functions from Sobolev spaces and quasi-optimal linear interpolators, Math. Inequal. Appl. 20 (2017) 441--457.
\bibitem {Fav} J. Favard, Sur l'interpolation, J. Math. Pures Appl. 19 (1940) 281--306.
\bibitem {FK-88} A. Fr\"olicher, A. Kriegl, Linear spaces and differentiation theory, Pure and Applied Mathematics, J. Wiley, Chichester, 1988.
\bibitem {JT} T. R. Jensen, B. Toft, Graph coloring problems, Wiley-Interscience Series in Discrete Mathe\-matics and Optimization. A Wiley-Interscience Publication. John Wiley $\&$ Sons, Inc., New York, 1995.
\bibitem {M} V. Maz’ya, Sobolev Spaces with Applications to Elliptic Partial Differential Equations, Springer,
Heidelberg, 2011.
\bibitem {Mer} J. Merrien, Prolongateurs de fonctions differentiables d'une variable r\'eelle, J. Math. Pures et Appl.  45 (1966) 291--309.
\bibitem {Sh-LMP-2018} P. Shvartsman, Extension criteria for homogeneous Sobolev spaces of functions of one variable, Rev. Mat. Iberoamericana (in press)
\bibitem {Sh-LV-2018} P. Shvartsman, Sobolev functions on closed subsets of the real line: long version, arXiv:1808.01467
\bibitem {Sh2} P. Shvartsman, Sobolev $W^1_p$-spaces on closed subsets of $\RN$, Adv. Math. 220 (2009) 1842--1922.
\bibitem {Sh5} P. Shvartsman, Whitney-type extension theorems for jets generated by Sobolev functions,  Adv. Math. 313 (2017) 379--469.
\bibitem {W2} H. Whitney, Differentiable functions defined in closed sets. I., Trans. Amer. Math. Soc. 36 (1934) 369--387.}
\end{thebibliography}
\end{document}